\documentclass[11pt,final]{iopart}

\expandafter\let\csname equation*\endcsname\relax
\expandafter\let\csname endequation*\endcsname\relax
\usepackage{amsmath}

\usepackage{algorithm,algorithmic}

\usepackage{amsthm}
\usepackage{amssymb}
\usepackage{graphicx}
\usepackage{fullpage}
\usepackage{xcolor}
\usepackage{color}
\usepackage{caption}
\usepackage{subcaption}
\usepackage{listings}
\usepackage{color}
\usepackage{hyperref}
\usepackage{tikz}
\usepackage{graphicx,scalerel}
\usepackage{placeins}

\DeclareMathOperator*{\argmin}{arg\,min}

\DeclareSymbolFontAlphabet{\amsmathbb}{AMSb}%

\usepackage[mathscr]{euscript}
\usepackage{color,amsmath,amssymb,mathrsfs}
\usepackage{graphicx}
\usepackage{algorithmic,algorithm}
\usepackage{tikz,pgfplots}
\usepackage{upgreek}
\usepackage{showkeys,cite}
\usepackage{multirow}
\usepackage{mathtools}
\usepackage[normalem]{ulem}
\usepackage{latexsym}
\DeclareMathAlphabet{\mathpzc}{OT1}{pzc}{m}{it}
\usepackage{url} 
\usepackage{tikz} 
\newcommand{\design}[2]{\bigg\{\!\! \begin{array}{c} {#1}\\{#2} \end{array}\!\!\!\bigg\}}

\definecolor{grassgreen}{RGB}{92,135,39}
\newcommand*{\transymb}{{\mkern-1.5mu\mathsf{T}}}
\renewcommand{\top}{\transymb}

\newcommand{\hilb}{\mathscr{H}}
\newcommand{\defeq}{\vcentcolon=}
\newcommand{\eqdef}{=\vcentcolon}

\renewcommand{\vec}[1]{{\mathchoice
                     {\mbox{\boldmath$\displaystyle{#1}$}}
                     {\mbox{\boldmath$\textstyle{#1}$}}
                     {\mbox{\boldmath$\scriptstyle{#1}$}}
                     {\mbox{\boldmath$\scriptscriptstyle{#1}$}}}}

\newcommand{\ran}{\mathsf{range}}

\newcommand{\trace}{\mathsf{tr}}
\newcommand{\eps}{\varepsilon}
\newcommand{\norm}[1]{\left\| {#1} \right\|}
\newcommand{\ip}[2]{{\left\langle {#1}, {#2} \right\rangle}}
\newcommand{\mip}[2]{\left\langle{#1}, {#2}\right\rangle_{\!\scriptscriptstyle{\text{M}}}}
\renewcommand{\mat}[1]{\mathbf{{#1}}}
\newcommand\restr[2]{{
  \left.\kern-\nulldelimiterspace 
  {#1}\vphantom{\big|} \right|_{#2}}}

\newcommand{\R}{\mathbb{R}}

\newcommand{\A}{\mathcal{A}}

\newcommand{\C}{\mathcal{C}}
\newcommand{\D}{\mathcal{D}}

\newcommand{\J}{\mathcal{J}}

\newcommand{\borel}{\mathscr{B}}

\newcommand{\var}{\mathbb{V}}
\newcommand{\GM}[2]{\mathcal{N}\!\left( {#1}, {#2}\right)}

\newcommand{\Cprior}{\mathcal{C}_{\text{pr}}}
\newcommand{\Cpost}{\mathcal{C}_{\text{post}}}
\newcommand{\ncov}{\mat{\Gamma}_{\!\text{noise}}}

\newcommand{\like}{\pi_{\text{like}}}
\newcommand{\obs}{\vec{y}}
\newcommand{\ipar}{m}
\newcommand{\iparpr}{m_{\text{pr}}}
\newcommand{\iparmap}{m_{\scriptscriptstyle\text{MAP}}}
\newcommand{\dpar}{\vec{m}}
\newcommand{\dparpr}{\dpar_{\text{pr}}}
\newcommand{\dparmap}{\vec{m}_{\scriptscriptstyle\text{MAP}}}
\newcommand{\priorcov}{\mat{\Gamma}_{\text{prior}} }
\newcommand{\postcov}{\mat{\Gamma}_{\text{post}} }
\newcommand{\priorm}{\mu_{\text{pr}}}
\newcommand{\postm}{\mu_{\text{post}}^{{\obs}}}
\newcommand{\M}{\mat{M}}                       

\newcommand{\ff}{\vec{f}}
\newcommand{\FF}{\EuScript{F}}
\newcommand{\bFF}{\mat{F}}

\newcommand{\Ns}{{n_s}}

\newcommand{\Nd}{{n_{\text{d}}}}

\renewcommand{\H}{\mathcal{H}}

\newcommand{\eeta}{\vec{\eta}}
\newcommand{\CM}{\mathscr{E}}
\newcommand{\eip}[2]{\left\langle{#1}, {#2}\right\rangle_{\!{\R^q}}}
\newcommand{\cip}[2]{\left\langle{#1}, {#2}\right\rangle_{\!\CM}}
\newcommand{\W}{\mat{W}}

\newcommand{\commentout}[1]{\iffalse {#1} \fi}

\newcommand{\DKL}[2]{D_\text{kl}\left({#1} \| {#2}\right)}

\newcommand{\PhiA}{\Phi_{\!\text{A}}}
\newcommand{\PhiD}{\Phi_{\text{D}}}
\newcommand{\Phic}{\Phi_{\text{c}}}
\newcommand{\hatPsiBR}{\widehat{\Psi}_{\text{risk}}}
\newcommand{\PsiBR}{\Psi_{\text{risk}}}
\newcommand{\dPhiA}{\boldsymbol{\Phi}_{\!\text{A}}}
\newcommand{\dPhiD}{\boldsymbol{\Phi}_{\text{D}}}
\newcommand{\dPhic}{\boldsymbol{\Phi}_{\text{c}}}
\newcommand{\dPhi}{\boldsymbol{\Phi}}
\newcommand{\PsiAG}{\Psi_{\!\text{A}}^{\text{G}}}
\newcommand{\PsicG}{\Psi_{c}^{\text{G}}}

\newcommand{\PsiA}{\Psi_{\!\text{A}}}

\newcommand{\irow}[1]{
  [\begin{matrix} #1 \end{matrix}]%
}

\begin{document}
\title{Optimal Experimental Design for Infinite-dimensional 
Bayesian Inverse Problems Governed by PDEs: A Review}

\author{Alen Alexanderian}

\address{
Department of Mathematics, North Carolina State University, Raleigh, NC.
}
\ead{\url{alexanderian@ncsu.edu}}

\date{\today}

\begin{abstract}
We present a review of methods for optimal experimental design (OED) for
Bayesian inverse problems governed by partial differential equations with
infinite-dimensional parameters.  The focus is on problems where one seeks to
optimize the placement of measurement points, at which data are
collected, such that the uncertainty in the estimated parameters is minimized. 
We present the mathematical foundations of OED in this context and
survey the computational methods for the class of OED problems under study. We
also outline some directions for future research in this area.

\end{abstract}



\section{Introduction}
Mathematical models of complex physical and biological systems play a crucial
role in understanding real world phenomena and making predictions.  Examples
include models of weather systems, ocean circulation, ice-sheet dynamics,
porous media flow, or spread of infectious diseases. Models governing complex
systems typically include a large number of parameters that are needed for a
full model specification.  Typically, some model parameters are uncertain and need
to be estimated using indirect measurements.  This is done by solving an
inverse problem~\cite{EnglHankeNeubauer96,Vogel02,Tarantola05} 
that uses the model and measurement 
data to estimate the unknown parameters.
Optimal experimental design (OED)~\cite{Pazman86,AtkinsonDonev92,Pukelsheim93,Ucinski05} 
comprises a critical component of parameter
estimation; it provides a rigorous framework to guide acquisition of 
data, using limited resources,
to construct model parameters with minimized uncertainty.

\begin{figure}[ht]\centering
\includegraphics[width=0.8\textwidth]{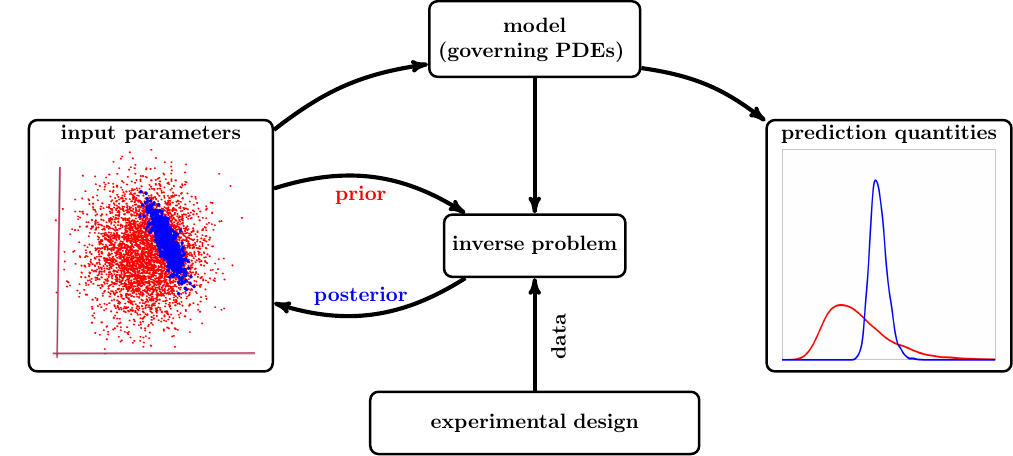}
\caption{An overview of Bayesian approach to inverse problems.
The Bayesian approach also enables making predictions with quantified and reduced 
uncertainties.}
\label{fig:diag_main}
\end{figure}

Models of complex systems are often described by systems of partial
differential equations (PDEs). Also, in many applications, the unknown
parameters to be estimated are functions, e.g., 
coefficient functions, initial states, or source
terms.  This review article is about OED for inverse problems governed by PDEs
with infinite-dimensional inversion parameters.  
We
focus on the Bayesian approach to inverse problems~\cite{Tarantola05,KaipioSomersalo06,Stuart10}, 
which provides a comprehensive framework for quantification of uncertainties inherent to
parameters and experimental data.  In this approach, we use measurement data
and a mathematical model to update the prior knowledge about unknown model
parameters. The prior knowledge is encoded in a prior distribution law for the
unknown parameters. The solution of a Bayesian inverse problem, known as the
posterior distribution, is a distribution law of the unknown parameters which
is consistent with measurement data, the model, and the prior distribution; see 
Figure~\ref{fig:diag_main}.
The study of the Bayesian approach to infinite-dimensional inverse problems
traces back to the article~\cite{Franklin70}.  Early follow up works include
\cite{Mandelbaum84}, \cite{PrenterVogel85},
\cite{LehtinenPaivarintaSomersalo89}, and \cite{Fitzpatrick91}. 
Significant  
progress has been made in this area in recent years; the
references~\cite{LassasSiltanen04,LassasSaksmanSiltanen09,Stuart10,
AgapiouLarssonStuart13,HelinBurger15,PinskiSimpsonStuartEtAl15,DashtiStuart17,Sprungk20} 
provide a small sample of literature from the last couple of decades.  There has also been
notable strides in computational methods for infinite-dimensional Bayesian
inverse problems; see, e.g.,~\cite{CotterRobertsStuartEtAl12,
HairerStuartVollmer14,PetraMartinStadlerEtAl14,IsaacPetraStadlerEtAl15,
CuiLawMarzouk16,BeskosGirolamiLanEtAl17,RudolfSprungk18,SaibabaBardsleyBrownAlexanderian19,BardsleyCuiMarzouk20}.

The OED problem seeks an experimental design, 
guiding data acquisition, that minimizes the uncertainty
in the estimated parameters.  Making this precise requires a discussion
of design criteria~\cite{ChalonerVerdinelli95,Ucinski05}.  There are various
such criteria, which attempt to measure the uncertainty in the estimated
parameters in different ways. For example, the so called A-optimal criterion
quantifies the average variance of the estimated parameters.

Measurement data, needed in parameter estimation, can be costly
or time consuming to acquire.  This can put severe limits on the amount of
experimental data that can be collected.  In such cases, a naive or otherwise
suboptimal experimental design entails waste of experimental resources 
and computing budget used to process the data, as well as 
possibly inaccurate parameter estimates.  Even in cases when data acquisition is not very
expensive, it is important to identify an optimal set of experiments. A
suboptimal data acquisition strategy could lead to processing datasets with
redundancies or, worse yet, one might miss important measurements leading to poor
parameter estimation. In short, OED is a crucial aspect of successful parameter
estimation and uncertainty quantification.

In Bayesian inversion and OED for infinite-dimensional inverse problems, it is
important to consider the problem in a function space setting, before
proceeding to discretization.  Not only does this lead to rigorous
formulations, it also is of practical importance when developing solution
methods: this guides the choice of prior measures that are meaningful for
infinite-dimensional parameters~\cite{Stuart10} and forces one to use
appropriate discretizations of the Bayesian inverse problem.  The latter is
important because if the discretization of the various components of a Bayesian
inverse problem is not performed correctly, the computed 
results will not converge to their infinite-dimensional
counterparts upon grid refinements; see~\cite{Bui-ThanhGhattasMartinEtAl13,
Bui-ThanhNguyen16}.
Considering the problem in an infinite-dimensional setting also  ensures
definition of OED criteria that are meaningful in infinite dimensions. These
criteria can then be discretized along with the  Bayesian inverse problem. This
way, the discretized OED criteria will have a meaningful infinite-dimensional
limit upon successive grid refinements.

OED for infinite-dimensional inverse problems is challenging from both
mathematical and computational points of view.  The definition and analysis of
OED criteria require tools from linear operator
theory~\cite{ReedSimon72,EinsiedlerWard17} and probability theory on Hilbert
spaces~\cite{PratoZabczyk92,Prato06}.  Computationally, one is faced with
optimization of expensive-to-evaluate (and differentiate) OED criteria.  Note
that the OED problem has the inverse problem as a sub-problem, which itself is
an extremely challenging problem.  These challenges are due to
high-dimensionality of the inversion parameters, upon discretization, and high
cost of simulating the governing model.  Addressing the computational
challenges of OED requires suitable approximations and computational methods
that maximally exploit the problem structure within the inverse problem. These
may include smoothing properties of the parameter-to-observable map or low-rank
structures in operators appearing in the definition of OED criteria.  

Broadly speaking, addressing mathematical and computational challenges of OED
for infinite-dimensional Bayesian inverse problems governed by PDEs requires an
interdisciplinary approach:  the solution methods involve an intricate blend of
methods and theories from mathematical and numerical analysis, inverse problem
theory, numerical methods for PDEs, probability theory in
function spaces, optimization, and uncertainty quantification. This article provides a review
of this framework, for the class of OED problems under consideration.
Our focus is primarily on situations where one seeks to optimize the locations
of points where measurement data are collected. A prime example 
is the problem of optimal sensor placement.

\textbf{A brief survey of literature}.
There is a rich body of literature devoted to OED for various classes of
parameter estimation problems.  Textbook references
include~\cite{Pazman86,AtkinsonDonev92,Pukelsheim93,Ucinski05,PronzatoPazman13,FedorovLeonov13}.
There are also a number of review articles on
OED~\cite{ChalonerVerdinelli95,Atkinson96,Clyde01,Muller05,RyanDrovandiMcGreeEtAl16}.  
Here we present a non-exhaustive sample of the literature on OED for models
governed by computationally intensive models.

Optimal design of experiments for inverse problems governed by ordinary
differential equations or differential-algebraic equations appear in numerous
works; examples include~\cite{BauerBockKorkelEtAl00,KorkelKostinaBockEtAl04,
BandaraSchloderEilsEtAl2009,ChungHaber12, BockKoerkelSchloeder13}. Recent works
also include OED for inverse problems constrained by PDEs; see
e.g.,~\cite{HoreshHaberTenorio10,AlexanderianPetraStadlerEtAl16,EtlingHerzog18,
HerzogRiedelUcinski18,Walter19,NeitzelPieperVexlerEtAl19,CastroDelosReyes20,Ucinski20}.
Bayesian approaches to OED for nonlinear inverse problems governed by
computationally challenging models were presented
in~\cite{HuanMarzouk13,HuanMarzouk14}; these articles use polynomial chaos
expansions and Monte Carlo to estimate the OED objective, and use a stochastic
optimization approach to compute the optimal design.  The use of polynomial
chaos surrogates makes this approach suitable for problems with low to moderate
parameter dimensions.  The
articles~\cite{LongScavinoTemponeEtAl13,LongMotamedTempone15,LongScavinoTemponeEtAl13}
use Laplace approximations to efficiently compute the expected information
gain, i.e., the D-optimal criterion, for nonlinear inverse problems.  The
article~\cite{WalshWildeyJakeman17} presents an approach for OED with an
alternate approach to statistical inverse problems in mind, known as the
consistent Bayes approach~\cite{ButlerJakemanWildey18}.  
The
articles~\cite{AlexanderianPetraStadlerEtAl14,AlexanderianPetraStadlerEtAl16,
AlexanderianGloorGhattas16,AlexanderianSaibaba18} concern OED for
infinite-dimensional Bayesian inverse problems. OED for infinite-dimensional 
problems is also treated in the PhD dissertation~\cite{Walter19}. 
See also the PhD dissertation~\cite{Herman20} that focuses on computational methods
for design of infinite-dimensional Bayesian linear inverse problems.
OED for large-scale problems has also been approached from a frequentist point
of view. The
articles~\cite{HaberHoreshTenorio08,HaberHoreshTenorio10,HoreshHaberTenorio10,HaberMagnantLuceroEtAl12}
provide excellent examples, where approaches based on Bayes risk minimization
are presented.  The article~\cite{ruthotto2017optimal} uses a similar
approach, but considers OED for inverse problems with state constraints.
The utility of OED has been explored in a wide range of application areas
involving large-scale inverse problems.  Examples include seismic waveform
inversion~\cite{DjikpesseKhodjaPrange12}, electrical impedance
tomography~\cite{HyvonenSeppanenStaboulis14}, borehole
tomography~\cite{HaberHoreshTenorio08}, contaminant
transport~\cite{AlexanderianPetraStadlerEtAl14,AlexanderianSaibaba18}
subsurface flow~\cite{AlexanderianPetraStadlerEtAl16}, natural gas
networks~\cite{YuZavalaAnitescu17},
thermomechanics~\cite{HerzogRiedelUcinski18}, iron loss in electrical
machines~\cite{HannukainenHyvonenPerkkio20}, and forensic
medicine~\cite{WeiserFreytagErdmannEtAl18}. 

\textbf{Article overview}. 
In Section~\ref{sec:motivation}, we present a few examples to motivate Bayesian
inversion and OED for models governed by differential equations with
infinite-dimensional parameters. Section~\ref{sec:HilbertBayes} outlines some
basics from probability theory and Bayesian inversion in Hilbert spaces. Theory
and methods for design of Bayesian linear inverse problems are reviewed in
Section~\ref{sec:oed_linear}. In that section, we also discuss a number of
tools and techniques that will be needed for nonlinear inverse problems as
well; these include definition of an experimental design, discretization
issues, randomized matrix methods, and sparsity control.  We discuss approaches
for design of Bayesian nonlinear inverse problems in
Section~\ref{sec:oed_nonlinear}. In our discussion of the numerical 
methods for OED, we also point to references containing extensive computational
results illustrating the effectiveness of the covered approaches. 
We conclude this review in
Section~\ref{sec:epilogue} with closing remarks and several 
directions for future work.

\section{Motivating Applications}\label{sec:motivation}
In this section, we present several motivating applications to provide 
intuition on design of infinite-dimensional inverse problems 
governed by PDEs. 

\subsection{Contaminant source identification}\label{sec:CSI}
In this example, we consider flow of a contaminant
in a geological formation. The present example is adapted
from~\cite{KovalAlexanderianStadler19}. Focusing on a horizontal 
cross-section of the medium, we consider a two 
dimensional domain $\D = [0,L_1] \times [0,L_2]$. 
The space-time evolution of the contaminant concentration, 
denoted by $u(\vec{x}, t)$, 
can be modeled
by a time-dependent advection-diffusion equation:
\begin{equation}\label{eq:ad-diff}
\begin{aligned}
\frac{\partial u}{\partial t}-\kappa\Delta{u}+\nabla \cdot (\vec v\, u) &= 0 \hspace{3mm} &&\mbox{   in } \D \times (0,T), \\
u(\cdot,0) &= m \hspace{3mm} &&\mbox{   in } \D, \\
(-\kappa\nabla{u} + u \,\vec{v}) \cdot\vec n  &= 0 \hspace{3mm} &&\mbox{   in }
\Gamma_{\!l} \times (0,T), \\
\kappa\nabla{u}\cdot \vec n &= 0 \hspace{3mm} &&\mbox{   in }
\partial \D \setminus \Gamma_{\!l} \times (0,T).
\end{aligned}\end{equation}
In the above problem, $\kappa > 0$ is the diffusion 
coefficient, $\vec{v}$ is the
velocity field, $T$ is the final
time,
and $m(\vec{x})$ is the initial concentration field. 
We have denoted the left edge of the domain by
$\Gamma_{\! l}$ and assume this part of
$\partial \D$ is  
impermeable, as modeled by the zero total flux condition.  The
homogeneous pure Neumann condition on the rest of the boundary
allows for advective flux. See~\cite{KovalAlexanderianStadler19}, for more details
on this model problem. 

In the present example, sensor measurements of the concentration are used to
estimate the initial concentration field $m(\vec{x})$.  This is an
example of a linear inverse problem, because the 
measurement data is a linear function of the parameter $m$. The OED
problem here seeks to specify the sensor locations so as to optimize the
statistical quality of the estimated initial concentration field.

\subsection{Permeability inversion in porous medium flow}\label{sec:pm_flow}
We consider the problem of estimating the permeability field of a
porous subsurface environment. We assume that a tracer substance flows
through the domain. Again, we focus on a two-dimensional 
geometry and let domain $\D$ be the unit square.
We denote by 
$\Gamma_0$, $\Gamma_1$, $\Gamma_2$, and $\Gamma_3$ 
the 
bottom, right, top, and left edges of the domain, respectively. 

The following equations model the flow of the tracer
through a medium that is saturated with a fluid: 
\begin{subequations}\label{equ:pde}
\begin{align} 
-\nabla\cdot(e^m \nabla p) &= 0 \quad 
\text{ in } \D \label{equ:pressure},\\
\frac{\partial u}{\partial t} - \nabla\cdot\big( \kappa \nabla u \big) + \nabla\cdot\big(
\vec{v}\, u\big)
&= g \quad \text{ in } [0,T]\times\D \label{equ:conc},\\
p &= p_1 \quad \text{on } \Gamma_1, \\
p &= p_2 \quad \text{on } \Gamma_3, \\
e^m \nabla p \cdot \vec{n} &= 0 \quad \text{ on } \Gamma_0 \cup \Gamma_2, \\
\nabla u \cdot \vec{n} &= 0 \quad \text{ on } [0,T]\times\{\Gamma_0\cup\Gamma_1\cup\Gamma_2\cup\Gamma_3\}, \\
u(0,\cdot) &= 0 \quad \text{ in } \D. 
\end{align} 
\end{subequations}
Here $p$ denotes the pressure field (of the saturating fluid), $m$
is the log-permeability field of the medium, 
$\vec{v} = -e^m\nabla p$ the Darcy velocity,
$u(t,\vec x)$ is the tracer concentration, 
$\kappa$ is the diffusion constant, and $g$ 
is a source term.
The Dirichlet pressure boundary conditions 
drive the flow. See~\cite{SunseriHartVanBloemenWaandersEtAl20}, 
from which this example is taken, for more details.

In this example, we use sensor measurements of 
pressure and concentration to estimate the log-permeability 
field $m(\vec{x})$. This is a nonlinear inverse problem; 
the function that maps the 
inversion parameter $m$ 
to the measurements is nonlinear. The OED problem here seeks to 
optimally place sensors that record concentration or pressure 
measurements. 

\subsection{Fault slip reconstruction in an earthquake}\label{sec:earthquake}
\begin{figure}[ht]\centering
\includegraphics[width=0.7\textwidth]{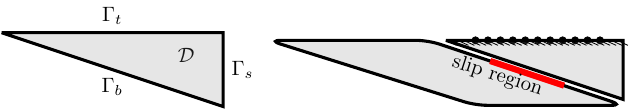}
\caption{Model geometry in fault slip reconstruction inverse problem.
Left: the computational domain. Right: idealized geometry depicting a subduction zone;
the GPS stations on the top boundary of the domain record displacement
during an earthquake.}
\label{fig:triangle}
\end{figure}

This problem, which is motivated by earthquake modeling,
concerns the reconstruction of the fault slip, during an earthquake.
The model problem here is adapted from~\cite[Chapter 2]{McCormack18}.
The 2D domain depicted in Figure~\ref{fig:triangle} is a 
cross-section of an idealized geometry modeling a subduction zone.
We represent the boundary of the domain $\D$
as $\partial \D = \Gamma_b\cup \Gamma_t\cup \Gamma_s$; see
Figure~\ref{fig:triangle}.
The interest here is on estimating the slip along $\Gamma_b$ based on point
measurements of elastic displacement on $\Gamma_t$.

Assuming a linear elastic equation for the
displacement $\vec u=\irow{u_1 & \! u_2}^\top$, we consider the governing model:
\[
-\nabla \cdot \vec \sigma(\vec u) = \vec 0, \quad \text{in } \D, 
\]
where 
$\sigma(\vec u) = 
\mu\,[\nabla \vec u + (\nabla \vec u)^\top] + \lambda\nabla \cdot \vec u \mat I$,
with $\mu$ and $\lambda$ denoting the L\'ame moduli.
The boundary conditions are:
\begin{align*}
\vec \sigma (\vec u)\vec n &= 0 \quad \text{ on } \Gamma_t,\\
\vec u  &= 0 \quad \text{ on } \Gamma_s,\\
\vec u\cdot \vec n &= 0 \quad \text{ on } \Gamma_b,\\
\delta \mat T (\vec \sigma(\vec u)\vec n)  + \mat T \vec u &= \vec m \quad \!\text{on } \Gamma_b.
\end{align*}
Here $\vec m$ is the inversion parameter, $\delta > 0$, and 
$\mat T$ is the tangential operator that
extracts the tangential components of a vector; i.e., 
$\mat T \vec u \defeq (\mat I -
\vec n \otimes \vec n)\vec u = \vec u - (\vec n^\top \vec u)\vec n$,
where $\vec n$ is the 
outward pointing unit normal on the boundary.

In this example, the inverse problem seeks to estimate 
$\vec{m}$ using displacement measurements at the top boundary. This is 
another example of a linear inverse problem. The OED problem here 
aims to optimally place the GPS stations that take displacement measurements 
on $\Gamma_t$.

\newcommand{\sign}{\phantom{-}}
\subsection{Parameter inversion in an epidemic model}\label{sec:epi}
While our focus is on inverse problems governed by PDEs, we also
present an example where the governing model is a system 
of ordinary differential equations.
Specifically, 
we consider a compartment model of the spread of a disease
known as the 
SEIRD model; see e.g.,~\cite{BrauerCastilloChavezFeng19} for 
more details on such models. This model 
tracks the populations of susceptible, exposed, infected, 
recovered, and dead individuals. The governing model 
can be written as 
\[
\begin{aligned}
S' &= -\beta S I / N, \\
E' &= \sign\beta S I / N - \sigma E,\\
I' &= \sign\sigma E - (\gamma + \delta) I,\\
R' &= \sign\gamma I,\\
D' &=\sign\delta I.
\end{aligned}
\]
Here $N = S+E+I+R+D$ and $\beta$, $\sigma$, $\gamma$, and $\delta$ are, respectively, the rates of
disease transmission, progression from exposed to infected, recovery, and
disease-related mortality. 
The parameters $\sigma$, $\gamma$, and $\delta$ are
scalars that depend on the specific disease.  However, $\beta$ is influenced by
the contact rate between individuals of a given population, which can change
over time. This can be, for example, as a result of social distancing policies
that reduce contact rate among individuals. Therefore, in general $\beta$ is
a function of time.  The inversion parameters here include the
infinite-dimensional parameter $\beta(t)$ and the scalar parameters $\sigma$,
$\gamma$, and $\delta$.

The inverse problem here seeks to estimate the unknown parameters using 
observations of infected people at a number of observation times.
An OED problem aims at finding optimal observation times.
Of course, in this inverse
problem one typically uses whatever data becomes available.  In
this case, OED is useful as it informs which measurement
times are important for reliable parameter estimation. If the available
datasets are missing data from key observation times, a practitioner would know
a priori that parameter estimation is subject to large uncertainties.

\section{Infinite-dimensional Bayesian inverse problems}\label{sec:HilbertBayes}
In this section, we outline background materials regarding Bayesian inversion
in an infinite-dimensional Hilbert space setting.  We first discuss the
requisite concepts regarding linear operators
(Section~\ref{sec:LinearOperators}) and probability measures
(Section~\ref{sec:borelmeas})  on infinite-dimensional Hilbert spaces.  Then,
we sketch the Bayesian formulation of an inverse problem in a Hilbert space in
Section~\ref{sec:BayesInv}.  For more details on the theoretical foundations of
Bayesian inversion in infinite dimensions and surveys of the related
literature, we refer the readers to~\cite{Stuart10,DashtiStuart17}.  

\subsection{Positive self-adjoint trace-class operators}\label{sec:LinearOperators}
Let $\hilb$ be an infinite-dimensional separable real Hilbert space equipped
with inner product $\ip{\cdot\,}{\cdot}$ and the corresponding
induced norm $\|\cdot\| = \ip{\cdot\,}{\cdot}^{1/2}$. A bounded linear operator 
$\A:\hilb \to \hilb$ is called selfadjoint if $\ip{\A u}{v} = \ip{u}{\A v}$
for all $u, v \in \hilb$. We call $\A$ positive if $\ip{\A u}{u} \geq 0$ 
for all $u \in \hilb$ and strictly positive if $\ip{\A u}{u} > 0$ for 
all nonzero $u \in \hilb$.  
A positive selfadjoint bounded linear operator $\A$ on $\hilb$ is called 
of trace-class if 
\begin{equation}\label{equ:trace}
\trace(\A) \defeq \sum_{k=1}^\infty \ip{\A e_k}{e_k} < \infty,
\end{equation}
where $\{ e_k \}_{k=1}^\infty$ is an orthonormal basis of $\hilb$. As 
is well-known~\cite{ReedSimon72,Conway00}, the value of the summation 
in~\eqref{equ:trace} does not depend on the 
choice of the orthonormal basis. Note also that if $\A$ is positive selfadjoint
and of trace-class, then it has an orthonormal basis of eigenvectors
$\{ v_k \}_{k=1}^\infty$ with corresponding (real, non-negative) eigenvalues
$\{ \lambda_k \}_{k=1}^\infty$, and $\trace(\A) = \sum_{k=1}^\infty \ip{\A v_k}{v_k}
= \sum_{k=1}^\infty \lambda_k$.

\subsection{Probability measures on a Hilbert space}\label{sec:borelmeas}
Here we recall some basics regarding probability measures on Hilbert spaces.

\textbf{Borel probability measures}. 
We let $\borel(\hilb)$ denote the Borel
$\sigma$-algebra on $\hilb$. 
A Borel probability measure on $\hilb$ is a
probability measure on the measurable space $(\hilb, \borel(\hilb))$. 
Let $\mu$ be a Borel probability measure on $\hilb$ that has  
finite
first and second moments. That is, $\int_\hilb \| z \| \, \mu(dz) < \infty$ and
$\int_\hilb \| z \|^2 \, \mu(dz) < \infty$.
The mean $\bar m$ of $\mu$ is an element of 
$\hilb$ that satisfies
\[
\ip{\bar m}{a} = \int_\hilb \ip{z}{a}\, \mu(dz), \quad \text{for all } a \in \hilb.
\]
The covariance operator of $\mu$ is a bounded linear operator
$\C:\hilb \to \hilb$ that is characterized as follows: 
\[
    \ip{\C a}{b} = \int_\hilb \ip{a}{z - \bar m}\ip{b}{z - \bar m} \, \mu(dz), 
    \quad a, b \in \hilb.
\]
Note that 
$\C$ is positive and selfadjoint; it is also straightforward 
to show that $\C$ is of trace-class, as seen 
from the following standard argument~\cite{Prato06}.
Let $\{e_k \}_{k=1}^\infty$ be 
a complete orthonormal set in $\hilb$. We have
\begin{equation}\label{equ:second_moment}
\begin{aligned}
\trace(\C) = \sum_{k=1}^\infty \ip{\C e_k}{e_k} 
&= \sum_{k=1}^\infty \int_\hilb \ip{e_k}{z - \bar{m}}^2 \, \mu(dz)
\\
&= \int_\hilb \sum_{k=1}^\infty \ip{e_k}{z - \bar{m}}^2 \, \mu(dz)
= \int_\hilb \|z - \bar m \|^2 \,\mu(dz) < \infty,
\end{aligned}
\end{equation}
where the interchange of the summation and the integral is justified by the
Lebesgue Monotone Convergence Theorem and the final equality is due to 
Parseval's identity.

We also recall that 
a Borel probability measure on $\mu$ on 
$(\hilb, \borel(\hilb))$ is uniquely characterized by
its Fourier transform~\cite[Proposition 1.7]{Prato06}, which is  
defined by
\[
\hat\mu(\xi) \defeq \int_\hilb e^{i\langle{x},{\xi}\rangle}\, \mu(dx),
\quad \xi \in \hilb.
\]

\textbf{Gaussian measures}.
An important class of probability measures encountered in Bayesian analysis in
function spaces is that of Gaussian measures.  We recall that $\mu$ is
a Gaussian measure on $(\hilb, \borel(\hilb))$ if for every $x \in \hilb$ the
linear functional $\ell(z) = \ip{x}{z}$, considered as random variable 
$\ell : (\hilb, \borel(\hilb), \mu) \to (\R, \borel(\R))$, is a Gaussian random
variable~\cite{DaPratoZabczyk14}.  Gaussian measures on $\R^n$ that have
strictly positive covariance operators can be characterized by the associated
probability density function (PDF). This PDF is none but the Radon--Nikodym
derivative~\cite{Williams1991} of the measure with the respect to the Lebesgue measure on $\R^n$.
However, in the infinite-dimensional Hilbert space setting, there is no
analogue of Lebesgue measure~\cite{Prato06}. In this case, we can 
characterize a Gaussian measure using its Fourier transform.
Namely, given $\bar m \in \hilb$ and a positive selfadjoint trace-class operator 
$\C$ on $\hilb$, 
the Gaussian measure $\mu = \GM{\bar m}{\C}$ is
the unique probability measure with,
\[
   \hat \mu(\xi) = \exp\Big\{ i\langle \bar m, \xi\rangle - \frac12 \langle Q\xi, \xi\rangle\Big\}, \quad \xi \in \hilb.
\]
For further details regarding Gaussian measures
on Hilbert spaces, see e.g.,\cite{DaPratoZabczyk02,Prato06,DaPratoZabczyk14}. 
See also~\cite{Bogachev98}, for a comprehensive treatment  of Gaussian 
measures on Banach spaces.

\textbf{Kullback-Leibler divergence.}
The Kullback-Leibler (KL) divergence provides a ``measure of distance'' between
two probability measures~\cite{Kullback1951}.  Note that the KL divergence is
not a metric, because it is non-symmetric and does not satisfy the triangle
inequality.  However, given probability measures $\mu_1$ and $\mu_2$, the
KL divergence from $\mu_1$ to $\mu_2$ is non-negative and is zero if and only
if the two measures are the same.  As discussed in Section~\ref{sec:criteria}, in 
Bayesian OED, the KL divergence is used to define a notion of
information gain.

For probability
measures on $\R^n$ that admit densities (i.e., PDFs) with respect to the Lebesgue measure,
we can define the KL divergence in terms of the densities. 
That
is, if $\mu_1$ and $\mu_2$ are probability measures with 
densities 
$\pi_1$ and $\pi_2$, respectively, 
the KL divergence from $\mu_1$ to $\mu_2$, denoted by $\DKL{\mu_1}{\mu_2}$,
is given by
\[
\DKL{\mu_1}{\mu_2}
= \int_{\R^n} \log\Big\{\frac{\pi_1(\vec{x})}{\pi_2(\vec{x})}\Big\}\, \pi_1(\vec{x}) \, d\vec{x}.
\]
However, in an infinite-dimensional Hilbert space, where there is no analogue
of the Lebesgue measure, we need to work with an abstract definition of KL
divergence, which we explain next.
Let $\mu_1$ and $\mu_2$ be two Borel probability measures on $\hilb$ and suppose $\mu_1$ is
absolutely continuous with respect to $\mu_2$. The KL divergence from $\mu_1$ to
$\mu_2$ is defined as 
\[ \DKL{\mu_1}{\mu_2} =
\int \log \Big\{\frac{d \mu_1}{d \mu_2}\Big\} \, d\mu_1, 
\] 
where
$\frac{d\mu_1}{d\mu_2}$ is the Radon-Nikodym derivative of $\mu_1$ with respect to
$\mu_2$.  In the case that $\mu_1$ is not absolutely continuous with respect to
$\mu_2$, we have $\DKL{\mu_1}{\mu_2} = +\infty$.

\textbf{$\hilb$-valued Random variables}.
Let $(\Omega, \Sigma, \mathbb{P})$ be a probability space~\cite{Williams1991} and let
$m:(\Omega, \Sigma, \mathbb{P}) \to (\hilb, \borel(\hilb))$ be a random variable.
The law $\mu$ of $m$ is a probability measure on $(\hilb, \borel(\hilb))$ 
that satisfies,
\[
\mu(E) = \mathbb{P}(m \in E) = \mathbb{P}(\{ \omega \in \Omega : 
m(\omega) \in E\}), \quad \text{ for } E \in \borel(\hilb).
\]
In this article, we assume $\hilb = L^2(\D)$ where $\D$ a bounded spatial
domain that has a sufficiently regular boundary.  (We can also
consider the case where $\D$ is a time interval.)  The inner product on $\hilb$
is the standard $L^2$ inner product.  In this case, for each $\omega \in
\Omega$, $m(\cdot,\omega)$ is a real-valued function on $\D$.  We can also consider $m$ as
a function $m: \D \times \Omega \to \R$, where for each $\vec{x} \in \D$,
$m(\vec{x}, \cdot)$ is a (scalar-valued) random variable; that is, $m$ is a random
field. As discussed below, in a Bayesian inverse problem governed by PDEs with
a random field parameter $m$, we seek to find a posterior distribution law for
$m$. This posterior law is conditioned on measurement data and is consistent
with a prior measure.  Herein, we use a common abuse of notation and use the
same letter to denote the random field parameter and the values taken by the
parameter. That is, we use $m$ to denote the random field parameter as well as
its realizations in $\hilb$.

\subsection{Bayesian inversion in an infinite-dimensional Hilbert space} 
\label{sec:BayesInv}
We consider the problem of inferring the distribution law of an uncertain 
parameter $m$, which takes values in $\hilb$, using measurement data $\obs$ and a model
$\ff(\ipar)$ that relates $m$ to data: 
\begin{equation}\label{equ:model}
    \obs = \ff(\ipar) + \vec{\eta}. 
\end{equation}
Here $\vec{\eta}$ is a random vector that quantifies measurement noise
and is assumed to be
independent of $\ipar$.
Evaluating
$\ff(\ipar)$, which we refer to as the parameter-to-observable map, 
for a given $\ipar \in \hilb$ involves
simulating a computational model, e.g., the governing PDEs, and applying 
a measurement operator that extracts observations from 
the solution of the governing model.
As mentioned before, 
we consider
the case that $\hilb = L^2(\D)$, where $\D$ is a bounded spatial 
domain with a sufficiently regular boundary.
In practical applications, one typically has access to 
finitely many (spatial or temporal) measurements. This is the case, 
in particular, when measurement data are collected at a set of sensors. 
Therefore, in the present work, 
we consider finite-dimensional observations $\obs \in \R^d$. 
However, we mention that a Bayesian inverse problem with 
infinite-dimensional observations can be 
considered as well, in which case $\ff$ is a mapping from $\hilb$ to a function
space; see e.g.,~\cite{DashtiStuart17}.

\textbf{The data likelihood}.
We denote by $\like(\obs | \ipar)$ the likelihood PDF, 
which describes the distribution of experimental data $\obs$ for a given
$\ipar \in \hilb$. Considering the 
common case of an additive Gaussian noise model, the random vector 
$\eeta$ in~\eqref{equ:model} is distributed according to 
$\GM{\vec{0}}{\ncov}$,
with $\ncov\in \R^{d\times d}$ being the noise covariance matrix. 
This implies that $\obs | \ipar
\sim \GM{\ff(\ipar)}{\ncov}$, and therefore, 
\begin{equation}\label{equ:likelihood} 
\like(\obs | \ipar) \propto \exp\left\{ -\frac12 \big(\ff(\ipar) - \obs\big)^\top \ncov^{-1} \big(\ff(\ipar) - \obs\big)\right\}.
\end{equation}

\textbf{The prior distribution law}.
Defining a prior measure that is meaningful in a function space
setting is non-trivial. This has been discussed in a number
of works; see e.g.,~\cite{Stuart10,LassasSaksmanSiltanen09,DashtiHarrisStuart12,
DashtiStuart17}. 
Herein, we assume 
a Gaussian prior 
$\priorm=\GM{\iparpr}{\Cprior}$ for the inference parameter, which is a common 
approach in the infinite-dimensional setting.
We assume that $\Cprior$ is a \emph{strictly} positive operator (it is also 
selfadjoint and trace-class).
A convenient way of defining the covariance 
operator is defining $\Cprior$ as the inverse of a
differential operator; see e.g.,~\cite{Bui-ThanhGhattasMartinEtAl13,Stuart10,
AlexanderianPetraStadlerEtAl16}.
The subspace $\CM = \ran(\Cprior^{1/2})$ is called
the Cameron--Martin space~\cite{Prato06} of the Gaussian measure $\priorm$.
This Cameron--Martin space is a dense subspace of $\hilb$ that 
is endowed with the inner product
\[
   \cip{x}{y} \defeq \ip{\Cprior^{-1/2} x}{\Cprior^{-1/2} y}, \quad x, y \in \CM.
\]
The prior mean $\iparpr$ is assumed to be an element of $\CM$. 

Specifying the prior can be viewed from a modeling standpoint.  The mean of the
prior can be thought of as our best guess for the inference parameter $\ipar$
before solving the inverse problem.
And we model correlation lengths and the pointwise variance of
the parameter through the covariance operator.  In the case of a covariance
operator defined in terms of inverse of a differential operator, the Green’s
function of the differential operator describes the correlation
structure~\cite{Bui-ThanhGhattasMartinEtAl13}.

\textbf{The Bayes Formula}.
Solving a Bayesian inverse problem amounts to finding the \emph{posterior}
distribution law, which we denote by $\postm$. The measure $\postm$  
describes a distribution law of the parameter $\ipar$
that is conditioned on measurement data.
Bayes' formula combines the ingredients of the Bayesian inverse problem 
and describes 
the relation between the prior
measure, the data likelihood, and this posterior measure.
In 
the infinite-dimensional Hilbert
space settings, Bayes' formula is given by \cite{Stuart10},
\begin{equation}\label{equ:bayes_abstract}
   \frac{d\postm}{d\priorm} \propto \like(\obs | \ipar),
\end{equation}
where the left hand side is the Radon--Nikodym derivative of
$\postm$ with respect to $\priorm$.  
The precise conditions, on the parameter-to-observable map 
and the prior measure, that ensure the above formula holds
is discussed in detail in~\cite{Stuart10}. Recall that the 
parameter-to-obsevable map, which is defined in terms of the governing model, 
enters the formulation of the Bayesian inverse problem
through the data likelihood; see~\eqref{equ:likelihood}.

Note that 
in the finite-dimensional setting the abstract form of the Bayes' formula above
can be reduced to the familiar form of Bayes' formula in terms of
PDFs.  Specifically, working in
finite-dimensions, with $\priorm$ and $\postm$ that are absolutely continuous
with respect to the Lebesgue measure $\lambda$, the prior and posterior
measures admit Lebesgue densities $\pi_\text{pr}$ and $\pi_\text{post}$,
respectively. Then, we note
\[
    \pi_\text{post}(\ipar | \obs) = 
    \frac{d\postm}{d\lambda}(\ipar) = 
    \frac{d\postm}{d\priorm}(\ipar)
    \frac{d\priorm}{d\lambda}(\ipar)
    \propto 
    \like(\obs | \ipar) 
    \pi_\text{pr}(\ipar).
\]

\textbf{The maximum a posteriori probability (MAP) point}.
The MAP point (or MAP estimator), which we denote by 
$\iparmap$, is a point estimate of the 
inversion parameter $\ipar$ that minimizes 
\begin{equation}\label{equ:MAP_cost}
\J(\ipar) \defeq \frac 12 \big(\ff(\ipar) - \obs\big)^\top
\ncov^{-1}\big(\ff(\ipar) - \obs\big) +
\frac12 \cip{\ipar - \iparpr}{\ipar - \iparpr},
\end{equation}
over the Cameron--Martin space $\CM$.
Note that $\iparmap$ depends on $\obs$; i.e., $\iparmap = \iparmap^\obs$.
For notational convenience, we suppress this dependence on $\obs$ in 
most of what follows.
The MAP point has an intuitive interpretation in the finite-dimensional 
setting---it maximizes the 
posterior PDF. As detailed 
in~\cite{Stuart10,DashtiLawStuartEtAl13,DashtiStuart17}, 
in the infinite-dimensional setting,
one can define the MAP point $\iparmap$ as a point $\ipar$
that maximizes the posterior probability of balls of radius $\eps$
centered at $\ipar$, in the limit as $\eps\to 0$. 
The existence of solutions to the problem 
\begin{equation}
\label{equ:inner-opt}
\min_{\ipar \in \CM} \J(\ipar) 
\end{equation}
can be established using standard variational arguments; see~\cite{Stuart10}.
However, in general, one cannot ensure uniqueness of the solution to~\eqref{equ:inner-opt}.

We call an inverse problem a \emph{linear inverse problem} if  
the parameter-to-observable map $\ff$ is linear; specifically, 
if $\ff(\ipar) = \FF\ipar$, with $\FF:\hilb \to \R^d$ a bounded linear transformation.
It is well-known~\cite{Stuart10} that the solution of a 
Bayesian linear inverse problem with Gaussian prior and noise models is 
a Gaussian posterior measure $\postm = \GM{\iparmap}{\Cpost}$, where 
\[
\Cpost = (\FF^* \ncov^{-1} \FF + \Cprior^{-1})^{-1}
\quad \text{and} \quad
\iparmap = \Cpost(\FF^* \ncov^{-1} \obs + \Cprior^{-1}\iparmap).
\]
In this case, the posterior covariance operator $\Cpost$, which does not
depend on measurement data $\obs$, provides a convenient  means to define measures of (posterior) 
uncertainty in the estimated parameter; cf.~Section~\ref{sec:criteria}.

\section{Design of linear inverse problems}\label{sec:oed_linear}
In this section, we focus on the classical case of Bayesian linear inverse
problems with Gaussian prior and noise models. This is an important special
case that deals with inverse problems with linear (or linearized)
parameter-to-observable maps.  The assumptions of Gaussianity on the prior and
noise, while by no means universal, are common, especially in the
infinite-dimensional setting.  In this case, as discussed in
Section~\ref{sec:HilbertBayes}, we have a Gaussian posterior measure.  The
posterior covariance operator, which appears in many of the standard OED
criteria, will depend on the choice of the experimental design, but will not
depend on observations.  This makes the formulation of the OED problem
straightforward. Though, the problem is not necessarily easy to solve
numerically.

We begin by recalling some of the commonly used OED criteria in the
infinite-dimensional setting, in Section~\ref{sec:criteria}.  We then proceed
to describe a specific setup of an OED problem---optimal sensor placement---in
Section~\ref{sec:sensors}.  We shall use this setup in our discussion of the
computational methods for OED.  The rest of the Section
describes discretization issues (Section~\ref{sec:discretization}), the numerical
optimization problem for finding an 
OED (Section~\ref{sec:OED_problem}), computational methods for computing
OED criteria
(Section~\ref{sec:challenges_methods}), sparsity control
(Section~\ref{sec:sparsity}), and convexity of OED criteria
(Section~\ref{sec:convexity}). We also discuss the greedy approach for solving
OED problems in Section~\ref{sec:greedy}.

\subsection{Common optimality criteria}\label{sec:criteria}
Here we discuss three commonly used OED criteria in the 
infinite-dimensional Hilbert space setting. 

\textbf{Bayesian A-optimality}.
In the finite-dimensional setting, an A-optimal design is one that minimizes
the trace of the posterior \emph{covariance matrix}~\cite{Ucinski05}. This
amounts to minimizing the average variance of an inference parameter vector.
As noted in~\cite{AlexanderianPetraStadlerEtAl14,
AlexanderianPetraStadlerEtAl16}, this notion extends naturally to the
infinite-dimensional setting.  Let $\ipar$ be an $\hilb$-valued inference
parameter in a Bayesian linear inverse problem.
Consider $m$ as a random variable 
$m:(\Omega, \Sigma, \mathbb{P}) \to (\hilb, \borel(\hilb), \postm)$,  
where $(\Omega, \Sigma, \mathbb{P})$ is a probability space. 
The pointwise posterior 
variance of $\ipar$ satisfies 
\begin{align*}
    \int_\D \var\{ m(\vec{x}, \cdot) \}  \, d\vec{x}
    &= \int_\D \int_\Omega [\ipar(\vec{x}, \omega) - \iparmap(\vec{x})]^2 
\,\mathbb{P}(d\omega)\, d\vec{x}
\\ 
&= \int_\Omega \int_\D [\ipar(\vec{x}, \omega) - \iparmap(\vec{x})]^2 
d\vec{x}\,\mathbb{P}(d\omega) \\
&= \int_\Omega \|\ipar(\cdot, \omega) - \iparmap\|^2 
\,
\mathbb{P}(d\omega) 
= \int_\hilb \| \ipar - \iparmap\|^2 \, \mu(dm) = \trace(\Cpost);
\end{align*}
here the interchange of the integrals is justified by Tonelli's 
theorem and the final equality follows from~\eqref{equ:second_moment}.
Thus, we see that the average posterior variance is proportional 
to $\trace(\Cpost)$. Therefore,
as in the finite-dimensional setting, 
the A-optimal criterion is given by
\[
   \PhiA = \trace(\Cpost),
\]
and minimizing $\PhiA$ amounts
to minimizing the average posterior variance of the inversion parameter.

Bayesian A-optimality can be viewed also from a 
decision-theoretic point of view. It is known~\cite{ChalonerVerdinelli95,AlexanderianGloorGhattas16}
that
\[
\int_\hilb \int_{\R^d}
\norm{\ipar - \iparmap^\obs}^2  \like( \vec{y} | \ipar) \, d\obs \, \priorm(d\ipar) = 
\PhiA.
\]
The expression on the left is known as the Bayes risk of the MAP point (with 
respect to the $L^2$ loss function).

\textbf{Bayesian c-optimality}.
In finite dimensions a c-optimal design minimizes the posterior variance of 
a linear combination of the inversion parameters; in the function
space setting, we consider the posterior variance of a weighted average
$\int_\D \ipar(\vec x) c(\vec x) \, d\vec x$.  Thus,
a Bayesian c-optimal design is one that minimizes the posterior variance of 
a linear functional 
\[
    \ell(\ipar) = \ip{c}{\ipar}, 
\]
for a fixed $c \in \hilb$. 
Note that
$\int_\hilb \ip{c}{\ipar} \postm(d\ipar) = \ip{c}{\iparmap}$,
and so the variance of $\ell$ is given by
\[
\var\{\ell\} = \int_\hilb (\ip{c}{\ipar} - \ip{c}{\iparmap})^2 \, \postm(d\ipar)
= \int_\hilb \ip{c}{\ipar - \iparmap}^2 \, \postm(d\ipar) = \ip{\Cpost c}{c}.
\]
This gives rise to the c-optimal criterion, 
\[
    \Phic = \ip{\Cpost c}{c}.
\]

We give another useful case where c-optimality is of interest.
Consider a scalar-valued prediction quantity of interest $p(m)$, where
$p:\hilb \to \R$ is a twice continuously differentiable function. 
Suppose we want a design that minimizes the posterior variance of $p(\ipar)$.
This can be difficult for a nonlinear function, as computing variance of 
$p(\ipar)$ requires sampling over $\ipar$. One can instead use a linearization
of $p$. 
Let $m_0$ be a point in $\hilb$, possibly a suitable guess for the inversion 
parameter $\ipar$. Consider the first order Taylor expansion,
\[
    p(\ipar) = p(m_0) + \ip{p'(m_0)}{\ipar - m_0} + o(\| \ipar - m_0\|).
\]
Using the local linear approximation 
$p_\text{lin}(\ipar) = p(m_0) + \ip{p'(m_0)}{\ipar - m_0}$,
we can compute the posterior variance of $p_\text{lin}$:
\[
    \var\{ p_\text{lin} \} = \ip{\Cpost p'(m_0)}{p'(m_0)}.
\]
In this case, the c-optimal design, with $c = p'(m_0)$, can be considered 
a \emph{goal-oriented} design criterion, measuring the \emph{approximate} posterior 
variance of a prediction quantity of interest.

\textbf{Bayesian D-optimality}.
In the finite-dimensional setting, D-optimality admits an intuitive geometric
interpretation: a  D-optimal design minimizes the volume of the
uncertainty ellipsoid. Mathematically, this is formulated as an optimization
problem that seeks to minimize the determinant of the posterior covariance
operator.  This, however, is not meaningful in infinite dimensions, because the
posterior covariance operator is a trace-class linear operator with eigenvalues
that accumulate at zero. 
The classical D-optimal criterion can also be understood from a decision-theoretic
point of view: a Bayesian D-optimal design is one that maximizes the expected 
information gain~\cite{ChalonerVerdinelli95}. This point
of view provides a roadmap for deriving the infinite-dimensional analogue of 
the D-optimal criterion~\cite{AlexanderianGloorGhattas16}.

The expected information gain can be defined as follows:
\[
    \text{expected information gain} :=  
\int_\hilb \int_{\R^d} \DKL{\postm}{\priorm}\, \like(\obs | \ipar) d\obs \, \priorm(d\ipar),
\]
where
$\DKL{\postm}{\priorm}$ is the KL divergence from posterior to 
prior: 
\[
    \DKL{\postm}{\priorm} = \int_\hilb \log\left\{ \frac{d\postm}{d\priorm} \right\}\, d\postm.
\]
As detailed in~\cite{AlexanderianGloorGhattas16}, we have
\[
   \int_\hilb \int_{\R^d} \DKL{\postm}{\priorm}\, \like(\obs | \ipar) d\obs \, \priorm(d\ipar) 
   = \frac12 \log\det (I + \tilde{\mathcal{H}}_m),
\]
where $\tilde{\mathcal{H}}_m = 
\Cprior^{1/2} \FF^* \ncov^{-1} \FF \Cprior^{1/2}$.
A D-optimal design is one that maximizes the 
expected information gain or, equivalently, minimizes 
\[
   \PhiD = -\log\det (I + \tilde{\mathcal{H}}_m).
\]

\subsection{Sensor placement}\label{sec:sensors}
Thus far, we have not been specific about the definition of an experimental
design. This is problem dependent.  To illustrate, we consider a few
examples. In an inverse problem of identifying the source of a contaminant 
(cf.\ Section~\ref{sec:CSI}), the
experimental design specifies the placement of sensors that take measurements
of the contaminant concentration.  In the example in
Section~\ref{sec:earthquake}, the design specifies the placements of GPS
stations taking measurements of displacement.  In an inverse problem governed
by an epidemic model, as in Section~\ref{sec:epi}, the design corresponds to the
observation times in which the number of infected individuals are recorded.
Another example comes from tomography, where the design could be the choice of
angles an object is hit with an x-ray source.  One can also have a sensor
placement problem with different types of sensors; e.g., in permeability
inversion in a porous media flow problem (cf.\ Section~\ref{sec:pm_flow}), one
can have sensors that take pressure measurements and sensors that take
concentration measurements.  An experimental design in that context specifies
the sensor locations and types.  This is a more complicated case of a sensor
placement problem, which we shall not discuss herein, but presents an
interesting avenue for future investigations. 

One approach to design of experiments is to formulate the problem as that of
selecting an optimal subset of a set of admissible
experiments~\cite{Ucinski05,HaberHoreshTenorio08,HaberMagnantLuceroEtAl12,
AlexanderianPetraStadlerEtAl14}.  Here we focus specifically on sensor
placement and describe a common approach to defining experimental designs in
this context. Further remarks on more general ways of defining experimental
designs will be discussed at the end of this subsection. 

We begin by fixing a set of points $\vec{x}_i$, $i = 1, \ldots, \Ns$, where
sensors can be placed.  These are the so called \emph{candidate sensor
locations}.  We assign a non-negative weight $w_i$ to each candidate location
$\vec{x}_i$ that, roughly speaking, indicates the importance of $\vec{x}_i$.
Then, as explained further below, we formulate the OED problem as  an
optimization problem in terms of $\vec{w} = \irow{w_1 & \!w_2 & \!\ldots &
\!w_\Ns}^\top$.   

We can interpret the weight vector $\vec{w}$ in various ways.  If experiments
are repeatable, the design weights can be used to guide  the number of times
to perform each experiment to reduce the corresponding measurement noise; see,
e.g., \cite{HaberHoreshTenorio08}.  In the context of sensor placement,
typically we seek  weight vectors containing zeros and ones only, indicating
whether or not to place a sensor at each of the candidate locations.  However,
due to combinatorial complexity of finding binary optimal design vectors, a
common approach is to consider a relaxation of the problem with weights $w_i
\in [0, 1]$. One then uses a suitable penalty method to obtain sparse and
binary optimal design vectors; this is discussed in section~\ref{sec:sparsity}.

\textbf{OED criteria as functions of $\vec w$}.
The weight vector $\vec{w}$ enters the 
Bayesian inverse problem~\eqref{equ:bayes_abstract} 
through the data likelihood~\cite{AlexanderianPetraStadlerEtAl14},
\[
   \frac{d\postm}{d\priorm} \propto \like(\obs | \ipar; \vec{w}). 
\]
Assuming $\ncov = \sigma^2\mat{I}$,
The $\vec{w}$-weighted data-likelihood is given by~\cite{AlexanderianPetraStadlerEtAl14},
\begin{equation}\label{equ:weighted-likelihood}
\like(\obs | \ipar; \vec{w}) \propto \exp\Big\{-\frac{1}{2\sigma^2}
(\FF\ipar - \vec{d})^\top\W(\FF\ipar - \vec{d})\Big\},
\end{equation}
where $\W \in \R^{d \times d}$ is a diagonal matrix with weights on its
diagonal. Recall that $d$ is the dimension of the vector of measurement data.
For time-dependent problems, $d$ equals the product of the number of candidate
sensor locations and the number of measurement times.  In that case, $\W$ is a
block diagonal matrix with each block being a diagonal matrix with $\vec{w}$ on
its diagonal.  Here, for simplicity, we assume measurements are taken at only
one time. That is, either the governing PDE is stationary, or in the
time-dependent setting, we take measurements at a final time; in this case, $d
= \Ns$. Note also that we have considered 
the common case of uncorrelated observations and assumed
a constant noise level at the measurement points. Modifying the present 
setup to allow for varying noise levels, but uncorrelated observations, is 
straightforward. For discussions of OED in problems with correlated measurements, 
see e.g.,~\cite{Muller07,LiuChepuriFardadEtAl16,Ucinski20}.

In the present setting, the $\vec{w}$-dependent posterior covariance operator is given by 
\begin{equation}\label{equ:postw}
\Cpost(\vec{w}) = (\sigma^{-2}\FF^* \W\FF + \Cprior^{-1})^{-1}.
\end{equation}
The
OED criteria defined above are now functions of $\vec{w}$. Namely,
\[
   \begin{aligned}
   \PhiA(\vec{w}) &= \trace\big(\Cpost(\vec{w})\big), \\
   \Phic(\vec{w}) &= \ip{\Cpost(\vec{w})c}{c}, \\
   \PhiD(\vec{w}) &= -\log\det (I + \sigma^{-2}\Cprior^{1/2} \FF^* \W \FF \Cprior^{1/2}). 
   \end{aligned}
\]

\textbf{Further remarks on defining experimental designs.}
We shall present the challenges and methods of computing OEDs for infinite-dimensional 
inverse problems in the context of the sensor placement problem as described above: 
find 
an optimal subset from an existing array of candidate sensor locations.
However, it is important to note that this setup 
is by no means the only point of view.  
Generally, one can define a design
(sensor placement) as a discrete probability measure~\cite{Ucinski05}. 
That is, one can define a design, denoted generically by $\xi$, as 
\begin{equation}\label{equ:design_general}
\xi \defeq \design{\vec x_1, \ldots, \vec x_{\Ns}}{p_1, \ldots, p_{\Ns}},
\quad \text{with } p_i \geq 0 \text{ and } \sum_{i = 1}^{\Ns} p_i = 1.
\end{equation}
The support points $\vec x_i$, $i = 1, \ldots, {\Ns}$ indicate sensor locations and
$p_i$, $i = 1, \ldots, \Ns$ indicate sensor weights. Note that fixing the support
points and optimizing over $p_i$'s is similar to the approach taken in the
beginning of this section (except there we did not require the weights to sum to
one). More generally, one can formulate the OED problem as an optimization
problem on the space of Borel measures over a subdomain where sensors can
be placed~\cite{KieferWolfowitz59,Pazman86,Pukelsheim93,Ucinski05}.  The idea is then to find optimal designs of the
form~\eqref{equ:design_general}. Existence of discrete optimal designs
is known in the case of finite-dimensional 
parameters; see e.g.,~\cite{Ucinski05,NeitzelPieperVexlerEtAl19}. The PhD 
Thesis~\cite{Walter19} generalizes the framework 
of optimal sensor placement, as an optimization problems on the set of
Borel measures, to the case of inverse problems governed by PDEs 
with infinite-dimensional parameters.

\subsection{Discretization}\label{sec:discretization}
Our focus being on Bayesian inverse problems governed by PDEs, we consider a
finite-element based discretization of the Bayesian inverse problem. 
A finite-element approach is natural in the context of inverse 
problems governed by PDEs~\cite{Bui-ThanhGhattasMartinEtAl13,
PetraMartinStadlerEtAl14,
Bui-ThanhNguyen16,
HoangQuekSchwab19,
Walter19}.
However, other approaches to discretization such as finite-difference
approxations or spectral approximations can be considered as well; 
see e.g.,~\cite{BeskosPinskiSanz-SernaEtAl11,RobertsRosenthal01}.
The
presentation here is based on the developments
in~\cite{Bui-ThanhGhattasMartinEtAl13} on finite-element 
discretization of Bayesian inverse problems. 

we consider a finite-element discretization of $m$, 
\[
    m_h(\vec{x}) = \sum_{i = 1}^n m_i \phi_i(\vec{x}),
\]
with $\phi_i$, $i= 1, \ldots, n$ Lagrange nodal basis functions. 
In this case,
identifying $m_h$ with $\dpar = [\begin{matrix} m_1 & m_2 & \cdots &
m_n\end{matrix}]^\top$, the discretized inversion parameter is the vector
$\dpar$ of the finite-element coefficients. Note that 
the discretized parameter dimension $n$ can be huge, especially in 
three-dimensional geometries.

For 
$u$ and $v$ in the span of $\phi_i$, $i = 1, \ldots, n$,
\[
    \ip{u}{v} = \int_\D u(\vec x) v(\vec x) \, d\vec{x} = 
    \sum_{i=1}^n \sum_{j = 1}^n u_i v_i \int_\D \phi_i(\vec{x}) \phi_j(\vec{x})\, d\vec{x}
    = \vec{u}^\top \mat{M} \vec{v} \eqdef \mip{\vec{u}}{\vec{v}},
\]
where $\M$ is the finite element mass matrix, $M_{ij} = \int_\D \phi_i(\vec{x})
\phi_j(\vec{x})\, d\vec{x}$.
The discretized parameter space is $\R^n$ equipped with 
the discretized $L^2(\D)$ inner product $\mip{\cdot}{\cdot}$
and norm $\norm{\cdot}_\mat{M} = \mip{\cdot}{\cdot}^{1/2}$.
The discretized 
parameter-to-observable map is a linear  transformation $\bFF:(\mathbb{R}^n, \mip{\cdot}{\cdot}) \to 
(\R^d, \ip{\cdot}{\cdot}_{\R^d})$, where $\ip{\cdot}{\cdot}_{\R^d}$ 
denotes the Euclidean inner product on $\R^d$. The discretized 
prior measure $\GM{\dparpr}{\priorcov}$ is obtained
by discretizing the prior mean and covariance operator, and 
the discretized posterior measure is given by $\GM{\dparmap}{\postcov}$, with
\[
\postcov = \left(\bFF^*\vec{\Gamma}^{-1}_\text{noise}\bFF + \priorcov^{-1}\right)^{-1}
\quad \text{and}
\quad 
\dparmap = \postcov\left(\bFF^*\vec{\Gamma}^{-1}_\text{noise}\obs+\priorcov^{-1}\vec{m}_\text{pr}
\right).
\]
\begin{figure}[ht]\centering
\includegraphics[width=.35\textwidth]{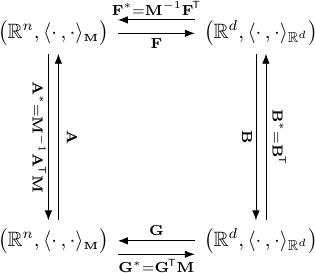}
\caption{Inner product spaces in a discretized Bayesian linear inverse problems 
and the adjoints of linear transformations between them.}
\label{fig:adjoints}
\end{figure}

The adjoints of the discretized linear operators involved in the Bayesian inverse problem need 
to also be defined carefully, as detailed in~\cite{Bui-ThanhGhattasMartinEtAl13}. 
Consider, for example,
a linear transformation $\mat{A}: (\R^n, \mip{\cdot\,}{\cdot}) \to (\R^n, \mip{\cdot\,}{\cdot})$.
The adjoint operator $\mat{A}^{\!*}$ satisfies, $\mip{\mat{A}\vec{u}}{\vec{v}}
= \mip{\vec{u}}{\mat{A}^{\!*}\vec{v}}$ for all $\vec{u}$ and $\vec{v}$ in $\R^n$.  We note that
\[
    \mip{\mat{A}\vec{u}}{\vec{v}}
    = (\mat{A}\vec{u})^\top \mat{M} \vec{v}
    = \vec{u}^\top \mat{A}^\top \mat{M} \vec{v}
    = \vec{u}^\top \mat{M} \mat{M}^{-1} \mat{A}^\top \mat{M} \vec{v}
    = \mip{\vec{u}}{\mat{M}^{-1} \mat{A}^\top \mat{M} \vec{v}},
\]
which shows that $\mat{A}^{\!*} = \mat{M}^{-1} \mat{A}^\top \mat{M}$. 
Note also that a selfadjoint operator $\mat{A}$ on $(\R^n,
\mip{\cdot\,}{\cdot})$ satisfies $\mat{A} = \mat{A}^{\!*}$.
For easy
reference, in Figure~\ref{fig:adjoints}, we summarize the definition of the
adjoint operators for linear transformations between spaces appearing in the
present setup of a Bayesian inverse problem. 
Notice that the prior and posterior covariance operators are selfadjoint with 
respect to the $L^2(\D)$ inner product, and their discretized versions 
will be selfadjoint with respect to the discretized $L^2(\D)$ inner product $\mip{\cdot\,}{\cdot}$.
See~\cite{Bui-ThanhGhattasMartinEtAl13},
for further details regarding the discretization of infinite-dimensional 
Bayesian inverse problems.

The OED criteria corresponding to the discretized Bayesian inverse problem are 
\begin{equation}\label{equ:OED_criteria}
   \begin{aligned}
   \dPhiA(\vec{w}) &= \trace\big(\postcov(\vec{w})\big), \\
   \dPhic(\vec{w}) &= \mip{\postcov(\vec{w})\vec{c}}{\vec{c}}, \\
   \dPhiD(\vec{w}) &= -\log\det (\mat{I} + \sigma^{-2}\priorcov^{1/2} \bFF^* \W \bFF \priorcov^{1/2}).
   \end{aligned}
\end{equation}

\subsection{The optimization problem for finding OEDs}\label{sec:OED_problem}
As described above, in our OED problem setup, a design is specified by 
a vector $\vec{w} \in [0, 1]^\Ns =: \mathcal{W}$.
The optimization problem for finding an optimal experimental design vector $\vec w$ can be 
formulated generically as 
\begin{equation}\label{equ:optim_problem}
\min_{\vec{w} \in \mathcal{W}} \dPhi(\vec{w}) + \gamma P(\vec{w}),
\end{equation}
where $\dPhi$ is a design criterion, e.g., ones listed in~\eqref{equ:OED_criteria}, 
$P(\vec{w})$ a penalty function used to promote sparsity and binary structure in $\vec{w}$, 
and $\gamma > 0$ a penalty parameter. 
This optimization problem can be solved with a gradient-based constrained optimization 
algorithm. 
Two important challenges associated with solving such optimization problems are (i)
the high cost of
evaluating the objective function and its gradient; and (ii) the need for an effective strategy in
choosing $P$ that leads to sparse and binary design vectors.

\subsection{Challenges in computing OED criteria and available approaches}\label{sec:challenges_methods}
To facilitate the discussion, 
let us consider A-optimality. The A-optimal criterion
$\dPhiA(\vec{w})$ is the trace of a high-dimensional operator---its dimension
is dictated by the discretized parameter dimension that can be 
in the thousands or hundreds of thousands in two- or three-dimensional 
computational domains.  
Applying this covariance operator to a vector is costly. 
Specifically, computing $\postcov(\vec{w}) \vec{v}$ requires computing,
\begin{equation}\label{equ:Gammapost_apply}
\left(\sigma^{-2}\bFF^*\W\bFF + \priorcov^{-1}\right)^{-1} \vec{v},
\end{equation}
which entails solving a high-dimensional linear system. Due to large-scale nature of the problem, 
matrix-free iterative methods such as the preconditioned conjugate gradient method should be used.
Note that computing the matrix-vector product (matvec) $\left(\sigma^{-2}\bFF^*\W\bFF + \priorcov^{-1}\right)\vec{v}$
requires a matvec with $\priorcov^{-1}$, as well as matvecs 
with $\bFF$ and $\bFF^*$. Computing a matvec with $\priorcov^{-1}$ can be challenging, 
but in problems where the prior covariance is defined in terms of an inverse of a 
differential operator, $\priorcov^{-1}$ is typically a sparse matrix that can be applied 
efficiently. In particular, no PDE solves are required for computing $\priorcov^{-1}\vec{v}$.
On the other hand, computing matvecs with $\bFF$ and $\bFF^*$, in inverse problems
governed by PDEs, require a forward and an adjoint PDE solve, respectively. 
Thus, the total cost of a matvec with $\postcov(\vec{w})$, in terms of PDE solves, will be 
two times the number of iterations of the iterative method used to 
compute~\eqref{equ:Gammapost_apply}.

Efficient means of applying the covariance operator $\postcov$, using low-rank
approximations in conjunction with the use of Sherman--Morrison--Woodbury
formula were presented in~\cite{Bui-ThanhGhattasMartinEtAl13}. The major cost
in this approach is that of computing a low-rank spectral decomposition of the
prior-preconditioned data-misfit Hessian $\sigma^{-2}\priorcov^{1/2}\bFF^*\W\bFF\priorcov^{1/2}$.  
Nevertheless, computing
$\trace(\postcov(\vec{w}))$ remains computationally challenging.  Also, this
needs to be repeated at every step of an optimization algorithm for finding an
OED.  Clearly building $\postcov(\vec{w})$ and computing its trace directly is
infeasible, due to the high computational cost. We will next
discuss methods for computing the OED objective.

\textbf{Methods based on Monte Carlo trace estimators.}
The articles~\cite{HaberHoreshTenorio08,HaberMagnantLuceroEtAl12}, propose the use of 
Monte Carlo trace estimators~\cite{Hutchinson90,AvronToledo11}. To recall briefly, 
a Monte Carlo estimator for $\trace(\mat{A})$, where $\mat{A}$ is a symmetric 
positive semidefinite matrix, is of the form
\begin{equation}\label{equ:MCtrace}
    \trace(\mat{A}) \approx \frac1N \sum_{i=1}^N \vec{v}_i^\top \mat{A} \vec{v}_i,
\end{equation}
where $\vec{v}_i$'s are realizations of a random vector with mean zero and identity 
covariance matrix. The reasoning behind these type of estimators stems
from the fact that 
$\mathbb{E}\left(\vec{v}^\top \mat{A} \vec{v}\right) = \trace(\mat{A})$, where
$\vec{v}$ is a random vector with mean zero and identity as the covariance matrix.
Examples include the Gaussian trace estimator, where entries of $\vec{v}$ are independent
and identically distributed (iid)
standard normal random variables, or the Hutchinson trace estimator~\cite{Hutchinson90}, where 
$v_i$ are iid Rademacher random variables.
(A Rademacher random variable takes values $\pm 1$ with probability of 
$1/2$ for each value.) It has been observed  
that a Monte Carlo estimator of the form~\eqref{equ:MCtrace} with a small $N$
can be effective in computing A-optimal OEDs. 

The idea of using a Monte Carlo trace estimator was also picked up in~\cite{AlexanderianPetraStadlerEtAl14},
for posterior covariance operators arising from discretization of an infinite-dimensional
Bayesian inverse problem. Namely, if $\vec{v}$ is a random $n$-vector 
with iid $\GM{0}{1}$ entries and with $\vec{z} = \M^{-1/2}\vec{v}$,
$\mathbb{E}\left(\mip{\vec{z}}{\postcov\vec{z}}\right) = \trace(\postcov)$; 
see~\cite[Proposition A.1.]{AlexanderianPetraStadlerEtAl14}. This leads
to an approximation of the form
\begin{equation}\label{equ:MCtraceM}
\trace\big(\postcov(\vec{w})\big) \approx \frac1N \sum_{i=1}^N \mip{\vec{z}_i}{\postcov(\vec{w})\vec{z}_i}, 
\end{equation}
with $\vec{z}_i = \M^{-1/2}\vec{v}_i$, where $\vec{v}_i$'s are realizations of
iid Gaussian random vectors. 

The use of Monte Carlo trace estimators provides a practical approach for
approximating the trace of the posterior covariance operator, but the
computational expense of the large number of matvecs with $\postcov$, incurred
over the course of the iterations of an optimization algorithm, can still be
formidable. This was addressed in~\cite{HaberMagnantLuceroEtAl12} by computing
a low-rank singular value decomposition of the forward operator $\bFF$.  Once a
low-rank approximation of $\bFF$ is available, $\postcov(\vec{w})$ applies can
be computed without any further forward/adjoint PDE solves.  A similar idea was
used in~\cite{AlexanderianPetraStadlerEtAl14}, where a low-rank SVD of the
\emph{prior-preconditioned} forward operator $\tilde\bFF = \bFF\priorcov^{1/2}$
was used. Due to the smoothing properties of the prior covariance operator,
$\tilde\bFF$ exhibits faster spectral decay, enabling further efficiency in
computing a low-rank approximation;
see~\cite{AlexanderianPetraStadlerEtAl14,AlexanderianSaibaba18}. 

The methods based on use of Monte Carlo trace estimators also facilitate
efficient computation of the gradient of the A-optimal criterion with respect
to $\vec{w}$~\cite{HaberMagnantLuceroEtAl12,AlexanderianPetraStadlerEtAl14}.
This enables an efficient optimization framework for computing A-optimal
designs for large-scale inverse problems.  It is also worth noting that the
approaches in~\cite{HaberMagnantLuceroEtAl12,AlexanderianPetraStadlerEtAl14}
can be used for computing a c-optimal criterion---$\dPhic$ is in the form of a
Monte Carlo estimator~\eqref{equ:MCtraceM} with $N = 1$.

The effectiveness of methods based on Monte Carlo trace-estimators for
computing A-optimal designs, for large-scale inverse problems, have been
demonstrated
in~\cite{HaberHoreshTenorio08,HaberMagnantLuceroEtAl12,AlexanderianPetraStadlerEtAl14}.
Specifically, this is demonstrated in a Borehole tomography problem
in~\cite{HaberHoreshTenorio08}, a magnetotelluric example problem and an
application to super-resolution reconstruction using MRI data
in~\cite{HaberMagnantLuceroEtAl12}, and a contaminat source identification
problem in~\cite{AlexanderianPetraStadlerEtAl14}.

\textbf{Randomized subspace iteration.}
Monte Carlo trace estimators are simple to implement and have been shown to
provide a practical way of tackling A-optimal design problems. However, they
exhibit slow convergence.  Significantly more accurate estimates of the trace can be
attained by exploiting a key problem structure---the often present low-rank
structure in the data misfit Hessian.  In~\cite{SaibabaAlexanderianIpsen17},
randomized trace estimators for trace and log-determinant were proposed that
rely on the concept of randomized subspace iteration.  These methods work well
for operators with large eigenvalue gaps or rapidly  decaying eigenvalues. 
Roughly speaking, the estimators based on subspace iteration work by
``projecting'' a matrix onto its dominant subspace. This process is sketched in
Algorithm~\ref{alg:randsvd}. 

\renewcommand{\algorithmicrequire}{\textbf{Input:}}
\renewcommand{\algorithmicensure}{\textbf{Output:}}

\begin{algorithm}[!ht]
\begin{algorithmic}[1]
\REQUIRE 
(i) Symmetric positive semi-definite matrix $\mat{A} \in\mathbb{R}^{n\times n}$ with
target rank $k$;
(ii) number of subspace iterations $q\geq 1$;
(iii) starting guess 
$\mat{\Omega} \in \mathbb{R}^{n\times \ell}$ with $k\leq \ell \leq n$, 
whose columns are random vectors with iid standard normal entries.
\ENSURE Matrix $\mat{T}\in\mathbb{R}^{\ell\times \ell}$.
\STATE Multiply $\mat{Y} = \mat{A}^q\mat{\Omega}$.
\STATE  Compute thin QR factorization $\mat{Y}=\mat{Q}\mat{R}$.
\STATE Compute $\mat{T} = \mat{Q}^\top\mat{A}\mat{Q}$.
\end{algorithmic}
\caption{Randomized subspace iteration (from~\cite{SaibabaAlexanderianIpsen17}).}
\label{alg:randsvd}
\end{algorithm}

With the output $\mat{T}$ of Algorithm~\ref{alg:randsvd}, we can approximate 
\[    
\trace(\mat{A}) \approx \trace(\mat{T}) \quad \text{and} \quad
\log\det(\mat{I} + \mat{A}) \approx
\log\det(\mat{I}+\mat{T}). 
\] 
Some comments on Algorithm~\ref{alg:randsvd} are in order.  In the first place,
the matrix $\mat{\Omega}$ can have entries from distributions other than
Gaussian; another possibility is to use $\mat{\Omega}$ with independent Rademacher
entries.  Moreover, in many cases, the choice of $q = 1$
is very effective in obtaining accurate estimates. 
For theoretical details of these estimators,
see~\cite{SaibabaAlexanderianIpsen17}.

The article~\cite{AlexanderianSaibaba18} presents methods for D-optimal design
of infinite-dimensional Bayesian linear inverse problems that are based on
randomized subspace iteration.  The work~\cite{HermanAlexanderianSaibaba19}
uses randomized subspace iteration for A-optimal design of linear inverse
problems.  Using these methods one can compute accurate approximations of the
OED objective and its gradient. Generally, as demonstrated
in~\cite{AlexanderianSaibaba18,HermanAlexanderianSaibaba19}, the methods based
on subspace iteration provide an excellent balance between accuracy,
computational efficiency, and ease of implementation.

\textbf{Trace and log-determinant evaluation in the measurement space}.
Recall that the discretized parameter dimension $n$ is typically very large.
In inverse problems where measurement data is collected at a set of sensors, 
the measurement dimension, in general, equals the number of candidate sensor locations
times the number of measurement times. In many cases, the discretized parameter 
dimension can be significantly larger than the measurement dimension. In such 
cases, it might be beneficial to reformulate the expressions for the OED criteria 
in such a way that the trace (or log-determinant) estimation is done in the measurement space, 
as we will illustrate next.

Let us consider the A-optimal criterion $\dPhiA$.  The weight-dependent 
posterior covariance operator can be written as
\begin{equation}\label{equ:posterior_cov_meas}
(\sigma^{-2}\bFF^* \W \bFF + \priorcov^{-1})^{-1}  = \priorcov - 
\sigma^{-2}\priorcov\bFF^*(\mat{I} + 
\sigma^{-2}\W \bFF \priorcov\bFF^*)^{-1}\W \bFF \priorcov;
\end{equation}
see~\cite{KovalAlexanderianStadler19}.
Since the prior covariance operator does not depend on $\vec{w}$, 
one can find an A-optimal design by minimizing the trace of the second
term in the right hand side of the above equation. After a simple manipulation 
we obtain the following OED objective:
\[
    \dPhi(\vec{w}) = \trace\left(
\sigma^{-2}(\mat{I} +
\sigma^{-2}\W \bFF \priorcov\bFF^*)^{-1}\W \bFF \priorcov^2 \bFF^*
\right).
\]
Note that the argument of the trace in the above expression is an operator 
on the measurement space. In cases where the dimension of the measurement 
space is significantly smaller than that of the discreized parameter space, 
this formulation can be useful for efficiently computing an A-optimal design.
This idea was used in~\cite{KovalAlexanderianStadler19} in context of 
A-optimal experimental design under uncertainty. Specifically, in~\cite{KovalAlexanderianStadler19}, A-optimal designs were computed in an inverse problem of estimating
an unknown initial state in a subsurface flow problem, under uncertainty in 
the flow field and the initial time.

A parallel development can be outlined in the case of D-optimality.
Consider $\dPhiD$ defined in~\eqref{equ:OED_criteria}. 
Using Sylvester's determinant identity, we can write,
\begin{equation}\label{equ:D-opt-new}
   \log\det (\mat{I} + \sigma^{-2}\priorcov^{1/2} \bFF^* \W \bFF \priorcov^{1/2})
   = \log\det(\mat{I} + \sigma^{-2} \W \bFF \priorcov \bFF^*).
\end{equation}
Therefore, again, we can consider formulating the OED problem with the 
expression on the right, which involves an operator 
defined on the measurement space, as the OED objective. 

\textbf{Adjoint-free approximate criteria}.
Another interesting problem structure revealed through these
``measurement space'' formulations is the 
possibility of eliminating the need for adjoint applies~\cite[Chapter 5]{Herman20}. 
For example, let us consider~\eqref{equ:D-opt-new} and focus on
the operator $\bFF \priorcov \bFF^*$ in the OED objective.
Suppose that the prior covariance operator has rapidly decaying eigenvalues 
allowing a low-rank approximation, $\priorcov \approx \mat{V}_r \mat{\Lambda}_r \mat{V}_r^*$.
Then, we can write 
\[
\bFF \priorcov \bFF^* \approx \bFF \mat{V}_r \mat{\Lambda}_r \mat{V}_r^* \bFF^*
= (\bFF \mat{V}_r\mat{\Lambda}_r^{1/2})(\bFF \mat{V}_r \mat{\Lambda}_r^{1/2})^*.
\]
In this case, one can consider building the operator $\tilde\bFF_r = \bFF \mat{V}_r \mat{\Lambda}_r^{1/2}$, 
at the cost of $r$ matvecs with $\bFF$ (forward applies), and using
the approximate D-optimal criterion:
\begin{equation}\label{equ:adjoint-free-Dopt}
   \widehat{\boldsymbol \Phi}_\text{D}(\vec{w}) =  -\log\det(\mat{I} + \sigma^{-2} \W \tilde\bFF_r \tilde\bFF_r^*).
\end{equation}
A key advantage of this formulation is that it does not require matvecs with $\bFF^*$. 
Note that in inverse problems governed by PDEs, computing the action of 
$\bFF^*$ requires an adjoint PDE solve. In applications where legacy solvers
are used, adjoint solvers might not be available and their implementation might not be feasible.
In such cases, the ``adjoint-free'' formulation~\eqref{equ:adjoint-free-Dopt} can be used.
A similar adjoint-free formulation can be derived for A-optimality, under the assumption 
that a low-rank approximation of $\priorcov$ is feasible. We refer to~\cite[Chapter 5]{Herman20}
for a detailed treatment, and numerical illustrations, of adjoint-free approaches
for computing A- and D-optimal designs.

\subsection{Sparsity control}\label{sec:sparsity}
We now return to the question of the choice of the penalty function $P(\vec w)$
in~\eqref{equ:optim_problem}.  In general, the choice of the penalty method
involves striking a balance between computational cost of the approach and the
ultimate goal of obtaining sparse and binary optimal design vectors.  A number
of techniques have been used to approximate the $\ell_0$-``norm'' to enforce
sparse and binary designs. For example, the authors 
of~\cite{HaberHoreshTenorio08} use an $\ell_1$-penalty
combined with a thresholding procedure.
In~\cite{AlexanderianPetraStadlerEtAl14}, a continuation approach is proposed
where a sequence of optimization problems with non-convex penalties that
successively approximate the $\ell_0$-``norm'' are solved.
In~\cite{HermanAlexanderianSaibaba19}, an approach based on reweighted
$\ell_1$-minimization is proposed. This also involves a continuation approach;
however, in each step, an optimization problem with a convex penalty is
solved. See also the related effort~\cite{YuZavalaAnitescu17}, in which a sum-up
rounding approach is proposed to obtain binary optimal designs. 

The OED problem in~\eqref{equ:optim_problem} can be solved 
via gradient-based optimization algorithms. Depending on the choice of 
the sparsity control approach, a number of optimization problems might 
need to be solved. For example, in~\cite{HermanAlexanderianSaibaba19}, where a 
reweighted $\ell_1$-minimization approach is used, one solves a sequence of optimization 
problems of the form, 
\[
    \min_{\vec w \in \mathcal{W}} \dPhi(\vec{w}) + \gamma \| \mat{D}_j \vec w\|_1,
    \quad j = 1, 2, \ldots.
\] 
Here $\mat{D}_j$'s are suitably chosen ``weighting matrices'';
see~\cite{HermanAlexanderianSaibaba19} for details. The efficient computational
methods for approximating the OED criteria and their gradients outlined above
can be used to accelerate the solution of such optimization problems, and
enable computing OEDs for infinite-dimensional Bayesian linear inverse
problems.

\subsection{Convexity of the common OED criteria}\label{sec:convexity}
Here we briefly comment on convexity properties of the OED criteria. 
We discuss the simpler case of c-optimality, but A- and D- optimality 
can be treated similarly; see e.g.,~\cite[Appendix B]{Ucinski05}. First we note
that the function 
\begin{equation}\label{equ:Gfun}
G(\mat{A}) := \mip{\mat{A}^{-1} \vec{c}}{\vec{c}}
\end{equation} 
is
strictly convex on the cone of strictly postive selfadjoint operators on 
$(\mat{R}^n, \mip{\cdot}{\cdot})$. This can be seen by the standard idea
of restricting $G$ to a line~\cite{BoydVandenberghe04}. Namely, consider
$g(t) = G(\mat{S} + t \mat{B})$, 
with $\mat{S}$ strictly positive and selfadjoint and $\mat{B}$ selfadjoint
on $(\R^n, \mip{\cdot}{\cdot})$; we consider $g(t)$ for values of $t$ such that 
$\mat{S} + t \mat{B}$ is strictly positve. Note that
\[
    0 < g(t) = \mip{(\mat{S} + t \mat{B})^{-1} \vec{c}}{\vec{c}} = 
           \mip{\mat{S}^{-1/2} (\mat{I} + t \mat{S}^{-1/2} \mat{B} \mat{S}^{-1/2})
                \mat{S}^{-1/2} \vec{c}}{\vec{c}}.
\]
We use the spectral decomposition of $\mat{S}^{-1/2} \mat{B} \mat{S}^{-1/2}$,
given by $\mat{S}^{-1/2} \mat{B} \mat{S}^{-1/2} = \sum_{i=1}^n \lambda_i
\vec{u}_i \otimes \vec{u}_i$ to write  $g(t) =  \sum_{i=1}^n (1+t\lambda_i)^{-1}
\mip{\vec{u}_i}{\mat{S}^{-1/2}\vec{c}}^2$.  This being a positive linear
combination of strictly convex functions, shows strict convexity of $g(t)$ and
subsequently, $G(\mat{A})$ with $\mat{A}$ in cone of strictly positive
selfadjoint operators on $(\R^n, \mip{\cdot}{\cdot})$.
 
Regarding convexity of the c-optimal objective, 
recall that $\postcov(\vec{w}) = \mat{H}(\vec{w})^{-1}$, with
$\mat{H}(\vec{w}) = \sigma^{-2}\bFF^*\mat{W} \bFF + \priorcov^{-1}$, 
and $\Phi(\vec{w}) = G(\mat{H}(\vec{w}))$. Here we consider 
$\vec{w} \in \R^\Ns_{\geq 0}$. Now, for $\vec{w}_1 \neq \vec{w}_2$,
we have, for $\alpha \in (0, 1)$, 
\begin{multline*}
    \dPhic\big(\alpha \vec{w}_1 + (1\!-\!\alpha) \vec{w}_2\big) 
    = G\big(\mat{H}(\alpha\vec{w}_1 + (1\!-\!\alpha)\vec{w}_2)\big) = 
    G\big(\alpha \mat{H}(\vec{w}_1) + (1\!-\!\alpha)\mat{H}(\vec{w}_2)\big)\\
     \leq 
    \alpha G\big(\mat{H}(\vec{w}_1)\big) + (1\!-\!\alpha)G\big(\mat{H}(\vec{w}_2)\big)  
     = 
    \alpha \dPhic(\vec{w}_1) + (1\!-\!\alpha) \dPhic(\vec{w}_2).
\end{multline*}
This shows convexity of $\dPhic$. Note that if $\mat{H}(\vec{w}_1) \neq \mat{H}(\vec{w}_2)$
the inequality in the penultimate step will be strict. Therefore, if we can
ensure $\mat{H}(\vec{w}_1) = \mat{H}(\vec{w}_2)$ implies $\vec{w}_1 =
\vec{w}_2$, we can conclude the strict convexity of $\dPhic(\vec{w})$.
This one-to-one property of $\mat{H}(\vec{w})$ can be obtained by putting
certain natural requirements on the forward operator, $\bFF \in \R^{\Ns
\times n}$.  Specifically, considering the typical case of $n > \Ns$, it
is straightforward to show that $\bFF$ having full row rank ensures that
$\mat{H}(\vec{w}_1) = \mat{H}(\vec{w}_2)$ imples $\vec{w}_1 =
\vec{w}_2$.

\subsection{Greedy approaches}\label{sec:greedy}
An alternative approach to sensor placements, which might be suitable
in some cases, is to follow a greedy procedure.  In this approach, we put
sensors sequentially, as outlined in Algorithm~\ref{alg:greedy}.

\begin{algorithm}[!ht]
\begin{algorithmic}[1]
\REQUIRE target number of sensors $K$. 
\ENSURE design vector $\vec w$.
\STATE set $\vec{w} = \vec{0}$,  $\mathcal{U} = \{1, \ldots, \Ns\}$, and
$\mathcal{S} = \emptyset$. 
\FOR{$k = 1$ \TO $K$} 
\STATE $\displaystyle i = \argmin_{j \in \mathcal{U} \setminus \mathcal{S}}
        \dPhi(\vec{w} + \vec{e}_j)$
\hfill
\COMMENT{$\vec{e}_j$ is the $j$th coordinate vector in $\R^\Ns$.}
\STATE $\mathcal{S} = \mathcal{S} \cup \{ i \}$.
\STATE $\vec{w} = \vec{w} + \vec{e}_i$
\ENDFOR
\end{algorithmic}
\caption{Greedy sensor placement.}
\label{alg:greedy}
\end{algorithm}
Theoretical justifications behind use of such an approach go back to
optimization of supermodular (or approximately supermodular) functions;
see~\cite{NemhauserWolseyFisher78,KrauseSinghGuestrin08,ChamonRibeiro17,ShulkindHoreshAvron18,
JagalurMarzouk20}. We recall that a function $f:2^U \to \R$, where $U$ is a finite set, 
is supermodular if 
\[
f(A \cup \{i\}) - f(S)
\leq f(B \cup \{i\}) - f(B) \quad \text{for all }
A \subset B \subset U, \text{and } i \in U\setminus B.
\]
The function $f$ is called submodular, if the inequality is reversed.

While the solution obtained from the greedy algorithm is sub-optimal, it can
provide good results in practice. The greedy algorithm is simple to implement, but
its cost, in terms of function evaluations, scales with the number $\Ns$ of
candidate sensor locations. Specifically, this requires $\mathcal{O}(K \Ns)$
function evaluations.  The efficient randomized methods for evaluating the OED
objective described above can be used to accelerate greedy sensor placement.

\section{Design of nonlinear inverse problems}\label{sec:oed_nonlinear}
In this section, we discuss OED for Bayesian nonlinear inverse problems.  That
is, we consider Bayesian inverse problems, where the parameter-to-observable
map is nonlinear.  For such problems, even with Gaussian prior and noise models, the posterior is
in general non-Gaussian; and, unlike the Gaussian linear inverse problems, no
closed-form expressions for measures of posterior uncertainty (design criteria)
are available.  A brute-force approach for computing measures of posterior
uncertainty would require generating samples from the posterior distribution
via an MCMC algorithm~\cite{Tarantola05}.  Such an approach to OED would be
infeasible for infinite-dimensional inverse problems governed by PDEs, because
an expensive MCMC procedure must be performed at every step of an optimization
algorithm.  This points to a fundamental challenge in design of large-scale
nonlinear inverse problems---the definition of suitable OED criteria whose
optimization is computationally tractable. 

A commonly used approach in classical works on design of nonlinear inverse
problems involves use of linearized models, leading to notions of locally 
optimum designs; see e.g.,~\cite{Atkinson96}. Such an approach can be done 
in a Bayesian setting as well, where a linearization point can be obtained 
based on prior information; e.g., one can use the prior mean for this.
Such an approach can also be done sequentially:  
alternate between estimating the uncertain parameters and obtaining an experimental 
design based on linearization at the current estimate; see e.g.,~\cite{KorkelBauerBockEtAl99}.

We discuss two approaches here that are developed in recent years 
to address OED for large-scale inverse problems: (i) compute an OED by minimizing the
(approximate) Bayes risk of the MAP point~\cite{HaberHoreshTenorio10}, and (ii) 
minimize approximate measure of posterior uncertainty obtained from 
a Laplace approximation to the posterior~\cite{AlexanderianPetraStadlerEtAl16}.
The former relies on ideas from decision theory to define the OED
objective.
Specifically, Bayes risk minimization
targets optimization of the statistical quality of the MAP point.
On the other hand, using a Laplace approximation enables 
incorporating approximate measures of posterior
uncertainty in the OED objective. These approaches, coupled with fast solvers, 
adjoint based gradient computation, and structure exploiting algorithms, can 
be turned into scalable algorithms for design of infinite-dimensional inverse 
problems~\cite{HaberHoreshTenorio10,AlexanderianPetraStadlerEtAl16,WuChenGhattas20}.

\subsection{Bayes risk minimization}
Given a design $\vec{w}$, we can compute the MAP point, $\iparmap(\vec{y}, \vec{w})$, 
by minimizing a functional of the form,
\begin{equation}\label{equ:MAP_cost_w}
\J_\vec{w}(\ipar, \obs) \defeq \frac {1}{2\sigma^2} \eip{\ff(\ipar) - \obs}{\W(\ff(\ipar) - \obs)} +
\frac12 \cip{\ipar - \iparpr}{\ipar - \iparpr},
\end{equation}
over the Cameron--Martin space $\CM$ (cf.\ section~\ref{sec:HilbertBayes}). 
Here
$\ff$ is the nonlinear parameter-to-observable map and $\W$ is the diagonal
matrix with entries of $\vec{w}$ on its diagonal. 
Note that we need data $\obs$ to find the MAP point, but we typically do not
have access to measurement data when solving the OED problem. In Bayes
risk minimization, one tackles this by considering an ``averaged 
criterion''.
The Bayes risk of the MAP 
point, with respect to the $L^2$ loss function, can be defined as 
\begin{equation}\label{equ:PsiBR}
    \PsiBR(\vec{w}) := \int_\hilb \int_{\R^d}  \| \iparmap(\obs, \vec{w}) - 
\ipar \|^2 \, \like(\obs | m)d\obs \, \priorm(dm).
\end{equation}
Numerically this criterion can be estimated via sample averaging. 
Namely, we can use draws $\{m_1, \ldots, m_\Nd\}$ from the prior distribution and 
we can take \emph{training} data samples 
\begin{equation}\label{equ:training_data}
\obs_i = \ff(m_i) + \vec\eta_i, \quad i = 1, \ldots, \Nd,
\end{equation}
where $\vec\eta_i$'s are draws from the noise distribution $\GM{\vec{0}}{\ncov}$.
The approximate Bayes risk then becomes
\begin{equation}\label{equ:hatPsiBR}
    \hatPsiBR(\vec{w}) := \frac{1}{\Nd} \sum_{i=1}^{\Nd} 
    \| \iparmap(\obs_i, \vec{w}) - \ipar_i \|^2. 
\end{equation}
Notice that each evaluation of this objective function requires solving for 
$\iparmap(\obs_i, \vec{w})$, $i = 1, \ldots, \Nd$. Thus, as formulated, the problem of minimizing
the Bayes risk is a bilevel optimization problem. In practice, 
this problem is formulated as~\cite{HaberHoreshTenorio10}
\begin{align}
&\min_{\vec{w}} \frac{1}{\Nd} \sum_{i=1}^{\Nd}
    \| \iparmap(\obs_i, \vec{w}) - \ipar_i \|^2 + \gamma P(\vec{w}),\label{equ:BR_obj}\\
&\text{where}\notag\\
&\J'_\vec{w}(\iparmap(\obs_i,\vec{w}), \obs_i) = 0, \quad i = 1, \ldots, \Nd.\label{equ:BR_const} 
\end{align}
Here $\J'_\vec{w}$ indicates the gradient of $\J_\vec{w}$ with respect to $m$,
and as before $P(\vec{w})$ is a sparsifying penalty function.  Note that  the
constraints~\eqref{equ:BR_const} are the first order optimality conditions for
the (inner) optimization problems of finding $\iparmap(\obs_i, \vec{w})$, $i =
1, \ldots, \Nd$. 

A few comments are in order. As one can see, the optimization problem~\eqref{equ:BR_obj} is 
a computationally challenging one: in particular, $\Nd$ inverse problems must be 
solved at each iteration of the optimization algorithm.  Therefore, in practice 
$\Nd$ cannot be very large. However, it is important to note that the inverse problem
solves can be performed in parallel. Computing the gradient of $\hatPsiBR(\vec{w})$ with
respect to $\vec{w}$ can be done efficiently using the adjoint method, making the cost
of gradient computation, in terms of the number of PDE solves, independent of discretized 
parameter dimension. 
Numerical illustrations of the effectiveness of this approach can be found
in~\cite{HaberHoreshTenorio10} in the context of an electromagnetic imaging
problem that uses direct current resistivity and magnetotelluric data.  See
also~\cite{ChungHaber12}, where Bayes risk minimization is used for design of
inverse problems governed by biological systems.

\subsection{OED criteria based on Laplace Approximation}
\textbf{A-optimal designs}. In this case, 
we seek designs that minimize the average posterior variance of the parameter
estimates  quantified by  the trace of the posterior covariance operator. 
To define an A-optimal criterion for nonlinear inverse problems we need to
overcome two fundamental challenges: (i) the posterior covariance operator
depends on the measurement data, which is not available a priori; and (ii) 
unlike the case of linear inverse problems, a
closed form expression for the posterior covariance operator is not available.

The first challenge can be addressed by considering an  averaged A-optimal
criterion; that is, we average over the set of all likely data, as done
in~\eqref{equ:PsiBR}. This results in
\begin{equation}\label{equ:oed-objective-general}
\PsiA(\vec{w}) := \int_\hilb \int_{\R^d} \trace(\Cpost(\vec{w}, \obs)) \, 
\like(\obs | m)d\obs \, \priorm(dm).
\end{equation}
Regarding the second challenge, while in principle
it is possible to estimate $\trace(\Cpost(\vec{w}, \obs))$ by generating samples from
the posterior distribution, this will be prohibitive for
large-scale problems. This calls for suitable approximations to the 
posterior distribution whose covariance operator admits a closed form 
expression. 

A commonly used tool, when working with large-scale nonlinear inverse problems is the 
Laplace approximation to the posterior. The Laplace approximation 
is a Gaussian approximation of the posterior measure,
with mean given by the MAP point
and covariance given by the inverse of the
Hessian (with respect to $\ipar$) operator $\H$ of $\J_\vec{w}$ 
in~\eqref{equ:MAP_cost_w}, evaluated at the MAP
point. More specifically, the Laplace approximation to the posterior is 
the Gaussian measure 
\begin{equation*}%
   \GM{\iparmap(\vec{w}, \obs)}{\H^{-1}\big(\iparmap(\vec{w}, \obs), \vec{w}, \obs\big)},
\end{equation*}
where $\H\big(\iparmap(\vec{w}, \obs), \vec{w}, \obs\big)$ is the Hessian
of~\eqref{equ:MAP_cost_w}.  In general, this Hessian depends on the design
vector $\vec{w}$ and data $\obs$ explicitly, as well as implicitly through the
MAP point~\cite{AlexanderianPetraStadlerEtAl16}.  Note that the Laplace
approximation to the posterior is exact when the parameter-to-observable map is
linear (and when we use Gaussian prior and noise models).  In a nonlinear
inverse problem, this approximation is suitable if the parameter-to-observable
map is well approximated by a linear approximation at the MAP point, over the
set of parameters with significant posterior probability.

Using this Gaussian approximation, we
can define an approximation $\PsiAG$ to the OED objective
defined in \eqref{equ:oed-objective-general}:
\begin{equation}\label{equ:oed-objective-Gaussian}
    \PsiAG(\vec{w}) = 
\int_\hilb \int_{\R^d}
\trace\left[\H^{-1}\big(\iparmap(\vec{w}, \obs), \vec{w}, \obs\big) \right]
\,\like(\obs | m)d\obs \, \priorm(dm).
\end{equation}  
In practice, this OED objective can be approximated via sample averaging, as
done in the case of Bayes-risk minimization; see~\eqref{equ:BR_obj}.  As in the
case of Bayes-risk minimization, finding designs that minimize (sample average
approximation to) $\PsiAG$ leads to a bilevel optimization problem.  An
additional challenge that needs to be addressed in optimization of $\PsiAG$ is
the need to estimate the trace of the inverse Hessian. This can be done
efficiently, for example, using randomized trace estimators.
See~\cite{AlexanderianPetraStadlerEtAl16} for a full elaboration of this
approach, where the effectiveness and scalability of this approach for OED is
demonstrated in a coefficient inversion problem motivated
by applications in porous medium flow. Below,  
we illustrate the optimization problem for finding an OED using the Laplace approximation
for the simpler case of Bayesian c-optimality.

\textbf{c-optimal designs}. 
In this case, we seek to optimize
posterior variance of a functional $\ip{c}{\ipar}$, where $c \in \hilb$. Analogously 
to~\eqref{equ:oed-objective-Gaussian}, we can define the Bayesian c-optimal criterion, 
in the nonlinear setting, as 
\[
\PsicG(\vec{w}) = \int_\hilb \int_{\R^d}  
\ip{\H^{-1}\big(\iparmap(\vec{w}, \obs), \vec{w}, \obs\big) c}{c} 
\,\like(\obs | m)d\obs \, \priorm(dm).
\]
We present the optimization problem for finding a c-optimal design
in an abstract form as follows, where also use sample averaging to approximate
$\PsicG$:
\begin{align}
&\min_{\vec{w}} \frac{1}{\Nd} \sum_{i=1}^{\Nd}
   \ip{z_i}{c}
+ \gamma P(\vec{w}),\label{equ:c-opt-obj}\\
&\text{where, for } i = 1, \ldots, \Nd\notag\\
&\J'_\vec{w}(\iparmap(\obs_i,\vec{w}), \obs_i) = 0, 
   \label{equ:c-opt-const_grad}
\\
&
\H\big(\iparmap(\vec{w}, \obs_i), \vec{w}, \obs_i\big) z_i = c.
\label{equ:c-opt-const_hess} 
\end{align}
As in the case of Bayes-risk minimization the training data vectors
$\obs_i$, $i = 1,\ldots, \Nd$ are generated according to~\eqref{equ:training_data}.
We point out that in PDE-based inverse problems~\eqref{equ:c-opt-const_grad}
will be replaced by the optimality system for the inverse problem, i.e., 
the PDEs describing the forward, adjoint, and gradient equations.
Additionally,~\eqref{equ:c-opt-const_hess}, which resembles a Newton step, 
is described by so called incremental state and adjoint equations and the equation 
describing the Hessian apply; see~\cite{AlexanderianPetraStadlerEtAl16}.

\textbf{Laplace approximations for Bayesian D-optimality}.
In Bayesian nonlinear inverse problems, computing the expected information
gain---the Bayesian D-optimal criterion---is challenging. The traditional
estimator for the expected information gain involves a double-loop Monte Carlo,
which usually requires large sample sizes; see~\cite{Ryan03,HuanMarzouk13}.
The Laplace approximation can be used for fast estimation of the expected
information gain. This idea has been exploited, for example,
in~\cite{LongScavinoTemponeEtAl13,LongMotamedTempone15,
BeckBenMansourEspathEtAl18}. These articles focus on inverse problems with
finite-dimensional parameters and seek to efficiently evaluate the D-optimal
objective. These efforts take fundamental steps towards development of a
scalable optimization framework for D-optimal design of infinite-dimensional
nonlinear Bayesian inverse problems. The recent work~\cite{WuChenGhattas20}
presents a computational framework for maximizing the expected
information gain, in infinite-dimensional Bayesian inverse problems, which 
utilizes Laplace approximations, low-rank approximations,
and greedy algorithms.

\section{Epilogue}\label{sec:epilogue}
OED for large-scale Bayesian inverse problems governed by PDEs is an exciting
and important area of research.  In our discussion, we highlighted common
problem formulations, challenges, approaches, and algorithms. The discussion,
while not exhaustive, reveals the richness of this field of research and also
points to a number of interesting directions for future work. We close our
discussion by listing a few such directions.

\textbf{OED for infinite-dimensional nonlinear inverse problems}.  OED for
nonlinear inverse problems governed by PDEs is challenging.  One can use
decision-theoretic criteria such as Bayes risk or expected information gain as
well as approximate measures of posterior uncertainty obtained via Laplace
approximations to the posterior.  There are other approximate criteria that
would be interesting to explore for large-scale problems. For example, the
review article~\cite{ChalonerVerdinelli95} lists a number of such approximate
measures.  From a computational mathematics point of view, with a given
computational budget, the choice of a specific approximate design criterion
must strike a balance in terms of computational complexity of the resulting OED
problem and the efficacy of the design criterion in measuring the statistical
quality of the estimated parameters. 

\textbf{OED under uncertainty}.
Typically, when solving large-scale inverse problems, one focuses on estimating
a specific set of model parameters and the remaining components of the model
are assumed known. This assumption is even more common when solving OED
problems. This, however, is not realistic in many cases.
Consider for example the contaminant source identification problem described in
Section~\ref{sec:CSI}, with the governing model~\eqref{eq:ad-diff}.  While the
inversion parameter there is the initial state, one can have (additional)
uncertainties in the diffusion coefficient, source term, velocity field, or
boundary conditions. Ignoring these additional uncertainties when designing
experiments can lead to vastly suboptimal designs in practice. A rational
approach for addressing this would be to consider uncertainties in the
additional model parameters in the OED problem formulation. This leads to the
formulation of the OED problem as an optimization under uncertainty problem. Initial
steps in this direction, for Bayesian linear inverse problems governed by PDEs,
are taken in~\cite{KovalAlexanderianStadler19,AlexanderianPetraStadlerEtAl20}. 
See also the recent work~\cite{FengMarzouk19}, where the authors consider 
\emph{focused} optimal experimental design for nonlinear models, 
where the goal is to maximize the 
expected information gain in targeted subsets of model parameters. 
Further developments in this
area present an interesting line of inquiry.

\textbf{Goal oriented OED}.
Design of an inverse problem should be performed with the ultimate goal of
solving the inverse problem in mind.  In some cases, solving an inverse problem
is merely an intermediate step in which a mathematical model is being
calibrated for the purposes of making predictions. In such cases, the
experimental design must be done with that final goal in mind. We refer to this
as  \emph{goal-oriented OED}.  This allows for a data collection strategy that
is tailored to the predictions. Not only does this allow for optimal use of
experimental resources, in many cases, the goal-oriented  design criteria can
be computationally easier to evaluate. The latter is due to the often
low-dimensionality of the prediction quantities of interest.  
Examples of such an approach
include~\cite{AttiaAlexanderianSaibaba18,HerzogRiedelUcinski18,Li19,ButlerJakemanWildey20}.  Further
work on goal-oriented design of inverse problems, especially for nonlinear
Bayesian inverse problems governed by complex physics systems, is an interesting
avenue of investigation.

\textbf{Switching or Mobile sensors}.
While we did not specifically discuss design of sensor networks for
time-dependent systems, the presented sensor placement setup can be extended to
a time-dependent setting also; see e.g.,~\cite{AlexanderianPetraStadlerEtAl14},
where an optimal placement of sensors taking measurements at a number of
observation times is considered.  In the time-dependent setting, it is also of
interest to consider switching or mobile sensors. In the case of switching (or
scanning) sensors, one seeks to find optimal activation protocols for sensors in
an already specified network of sensors: an optimal subset of sensors should be
activated at each observation time.  Motivations for this include reducing the
amounts of data to be processed and managing the cost of operating the
sensors~\cite{Ucinski05}.  In the case of mobile sensors one considers
measurement devices placed on monitoring cars or drones. A thorough treatment of
designing scanning or mobile sensor configurations is given in~\cite[Chapter 4]{Ucinski05}.  
Methods for design of such sensor configurations for
large-scale Bayesian inverse problems governed by time-dependent PDEs is an
interesting and important area for further research.

\section*{Acknowledgements}
I would like to thank Georg Stadler for numerous helpful discussions during
preparation of this article and for carefully reading through a draft of
this work.  Thanks also to No\'emi Petra for helpful discussions and feedback.  I
am grateful to Karina Koval for carefully reading through a draft of this
article and for help with proofreading.   

This work was supported in part by the National Science Foundation under 
the grant DMS-1745654.

\section*{References}
\bibliographystyle{abbrv}
\bibliography{refs}

\begin{thebibliography}{100}

\bibitem{AgapiouLarssonStuart13}
S.~Agapiou, S.~Larsson, and A.~M. Stuart.
\newblock Posterior contraction rates for the {B}ayesian approach to linear
  ill-posed inverse problems.
\newblock {\em Stochastic Processes and their Applications},
  123(10):3828--3860, 2013.

\bibitem{AlexanderianGloorGhattas16}
A.~Alexanderian, P.~J. Gloor, and O.~Ghattas.
\newblock On {B}ayesian {A}-and {D}-optimal experimental designs in infinite
  dimensions.
\newblock {\em Bayesian Analysis}, 11(3):671--695, 2016.

\bibitem{AlexanderianPetraStadlerEtAl14}
A.~Alexanderian, N.~Petra, G.~Stadler, and O.~Ghattas.
\newblock {A}-optimal design of experiments for infinite-dimensional {B}ayesian
  linear inverse problems with regularized $\ell_0$-sparsification.
\newblock {\em SIAM Journal on Scientific Computing}, 36(5):A2122--A2148, 2014.

\bibitem{AlexanderianPetraStadlerEtAl16}
A.~Alexanderian, N.~Petra, G.~Stadler, and O.~Ghattas.
\newblock A fast and scalable method for {A}-optimal design of experiments for
  infinite-dimensional {B}ayesian nonlinear inverse problems.
\newblock {\em SIAM Journal on Scientific Computing}, 38(1):A243--A272, 2016.

\bibitem{AlexanderianPetraStadlerEtAl20}
A.~Alexanderian, N.~Petra, G.~Stadler, and I.~Sunseri.
\newblock Optimal design of large-scale bayesian linear inverse problems under
  reducible model uncertainty: good to know what you don't know.
\newblock {\em SIAM/ASA Journal on uncertainty quantification}, Accepted, 2020.
\newblock \url{https://arxiv.org/abs/2006.11939}.

\bibitem{AlexanderianSaibaba18}
A.~{Alexanderian} and A.~{Saibaba}.
\newblock {Efficient D-optimal design of experiments for infinite-dimensional
  {B}ayesian linear inverse problems}.
\newblock {\em SIAM Journal on Scientific Computing}, 40(5), 2018.

\bibitem{Atkinson96}
A.~Atkinson.
\newblock The usefulness of optimum experimental designs.
\newblock {\em Journal of the Royal Statistical Society: Series B
  (Methodological)}, 58(1):59--76, 1996.

\bibitem{AtkinsonDonev92}
A.~C. Atkinson and A.~N. Donev.
\newblock {\em Optimum Experimental Designs}.
\newblock Oxford, 1992.

\bibitem{AttiaAlexanderianSaibaba18}
A.~{Attia}, A.~{Alexanderian}, and A.~K. {Saibaba}.
\newblock {Goal-Oriented Optimal Design of Experiments for Large-Scale
  {B}ayesian Linear Inverse Problems}.
\newblock {\em Inverse Problems}, 34(9), 2018.

\bibitem{AvronToledo11}
H.~Avron and S.~Toledo.
\newblock Randomized algorithms for estimating the trace of an implicit
  symmetric positive semi-definite matrix.
\newblock {\em Journal of the ACM (JACM)}, 58(2):17, April 2011.

\bibitem{BandaraSchloderEilsEtAl2009}
S.~Bandara, J.~P. Schl{\"o}der, R.~Eils, H.~G. Bock, and T.~Meyer.
\newblock Optimal experimental design for parameter estimation of a cell
  signaling model.
\newblock {\em PLoS computational biology}, 5(11), 2009.

\bibitem{BardsleyCuiMarzouk20}
J.~M. Bardsley, T.~Cui, Y.~M. Marzouk, and Z.~Wang.
\newblock Scalable optimization-based sampling on function space.
\newblock {\em SIAM Journal on Scientific Computing}, 42(2):A1317--A1347, 2020.

\bibitem{BauerBockKorkelEtAl00}
I.~Bauer, H.~G. Bock, S.~Körkel, and J.~P. Schlöder.
\newblock Numerical methods for optimum experimental design in {DAE} systems.
\newblock {\em Journal of Computational and Applied Mathematics}, 120(1):1--25,
  2000.

\bibitem{BeckBenMansourEspathEtAl18}
J.~Beck, B.~M. Dia, L.~F. Espath, Q.~Long, and R.~Tempone.
\newblock Fast {B}ayesian experimental design: Laplace-based importance
  sampling for the expected information gain.
\newblock {\em Computer Methods in Applied Mechanics and Engineering},
  334:523--553, 2018.

\bibitem{BeskosGirolamiLanEtAl17}
A.~Beskos, M.~Girolami, S.~Lan, P.~E. Farrell, and A.~M. Stuart.
\newblock Geometric {MCMC} for infinite-dimensional inverse problems.
\newblock {\em Journal of Computational Physics}, 335:327--351, 2017.

\bibitem{BeskosPinskiSanz-SernaEtAl11}
A.~Beskos, F.~J. Pinski, J.~M. Sanz-Serna, and A.~M. Stuart.
\newblock Hybrid {M}onte {C}arlo on {H}ilbert spaces.
\newblock {\em Stochastic Processes and their Applications}, 121:2201--2230,
  2011.

\bibitem{BockKoerkelSchloeder13}
H.~G. Bock, S.~K\"orkel, and J.~P. Schl\"oder.
\newblock Parameter estimation and optimum experimental design for differential
  equation models.
\newblock In H.~G. Bock, T.~Carraro, W.~J\"ager, S.~K\"orkel, R.~Rannacher, and
  J.~P. Schl\"oder, editors, {\em Model Based Parameter Estimation}, volume~4
  of {\em Contributions in Mathematical and Computational Sciences}, pages
  1--30. Springer Berlin Heidelberg, 2013.

\bibitem{Bogachev98}
V.~I. Bogachev.
\newblock {\em Gaussian {M}easures}.
\newblock Mathematical Surveys and Monographs. Americal Mathematical Society,
  1998.

\bibitem{BoydVandenberghe04}
S.~Boyd and L.~Vandenberghe.
\newblock {\em Convex optimization}.
\newblock Cambridge University Press, Cambridge, 2004.

\bibitem{BrauerCastilloChavezFeng19}
F.~Brauer, C.~Castillo-Chavez, and Z.~Feng.
\newblock {\em Mathematical Models in Epidemiology}.
\newblock Springer, 2019.

\bibitem{Bui-ThanhGhattasMartinEtAl13}
T.~Bui-Thanh, O.~Ghattas, J.~Martin, and G.~Stadler.
\newblock A computational framework for infinite-dimensional {B}ayesian inverse
  problems {P}art {I}: {T}he linearized case, with application to global
  seismic inversion.
\newblock {\em SIAM Journal on Scientific Computing}, 35(6):A2494--A2523, 2013.

\bibitem{Bui-ThanhNguyen16}
T.~Bui-Thanh and Q.~P. Nguyen.
\newblock {FEM-based discretization-invariant MCMC methods for PDE-constrained
  Bayesian inverse problems}.
\newblock {\em Inverse Problems \& Imaging}, 10(4):943, 2016.

\bibitem{ButlerJakemanWildey18}
T.~Butler, J.~Jakeman, and T.~Wildey.
\newblock Combining push-forward measures and {B}ayes' rule to construct
  consistent solutions to stochastic inverse problems.
\newblock {\em SIAM Journal on Scientific Computing}, 40(2):A984--A1011, 2018.

\bibitem{ButlerJakemanWildey20}
T.~Butler, J.~Jakeman, and T.~Wildey.
\newblock Optimal experimental design for prediction based on push-forward
  probability measures.
\newblock {\em Journal of Computational Physics}, 416, 2020.

\bibitem{CastroDelosReyes20}
P.~Castro and J.~C. De~los Reyes.
\newblock A bilevel learning approach for optimal observation placement in
  variational data assimilation.
\newblock {\em Inverse Problems}, 36(3):035020, 2020.

\bibitem{ChalonerVerdinelli95}
K.~Chaloner and I.~Verdinelli.
\newblock Bayesian experimental design: A review.
\newblock {\em Statistical Science}, 10(3):273--304, 1995.

\bibitem{ChamonRibeiro17}
L.~Chamon and A.~Ribeiro.
\newblock Approximate supermodularity bounds for experimental design.
\newblock In {\em Advances in Neural Information Processing Systems}, pages
  5403--5412, 2017.

\bibitem{ChungHaber12}
M.~Chung and E.~Haber.
\newblock Experimental design for biological systems.
\newblock {\em SIAM Journal on Control and Optimization}, 50(1):471--489, 2012.

\bibitem{Clyde01}
M.~A. Clyde.
\newblock Experimental design: A {B}ayesian perspective.
\newblock {\em International Encyclopia Social and Behavioral Sciences}, 2001.

\bibitem{Conway00}
J.~B. Conway.
\newblock {\em A course in operator theory}.
\newblock American Mathematical Soc., 2000.

\bibitem{CotterRobertsStuartEtAl12}
S.~Cotter, G.~Roberts, A.~Stuart, and D.White.
\newblock {MCMC} methods for functions: modifying old algorithms to make them
  faster.
\newblock {\em Statistical Science}, pages 424--446, 2013.

\bibitem{CuiLawMarzouk16}
T.~Cui, K.~J. Law, and Y.~M. Marzouk.
\newblock Dimension-independent likelihood-informed {MCMC}.
\newblock {\em Journal of Computational Physics}, 304:109--137, 2016.

\bibitem{Prato06}
G.~Da~Prato.
\newblock {\em An Introduction to Infinite-dimensional Analysis}.
\newblock Universitext. Springer, 2006.

\bibitem{DaPratoZabczyk02}
G.~Da~Prato and J.~Zabczyk.
\newblock {\em Second-order partial differential equations in Hilbert spaces}.
\newblock Cambridge University Press, 2002.

\bibitem{DaPratoZabczyk14}
G.~Da~Prato and J.~Zabczyk.
\newblock {\em Stochastic equations in infinite dimensions}, volume 152.
\newblock Cambridge university press, 2014.

\bibitem{DashtiHarrisStuart12}
M.~Dashti, S.~Harris, and A.~Stuart.
\newblock Besov priors for {B}ayesian inverse problems.
\newblock {\em Inverse Problems and Imaging}, 6(2):183--200, 2012.

\bibitem{DashtiLawStuartEtAl13}
M.~Dashti, K.~J. Law, A.~M. Stuart, and J.~Voss.
\newblock {MAP} estimators and their consistency in {B}ayesian nonparametric
  inverse problems.
\newblock {\em Inverse Problems}, 29, 2013.

\bibitem{DashtiStuart17}
M.~Dashti and A.~M. Stuart.
\newblock The {B}ayesian approach to inverse problems.
\newblock In R.~Ghanem, D.~Higdon, and H.~Owhadi, editors, {\em Handbook of
  Uncertainty Quantification}. Spinger, 2017.

\bibitem{DjikpesseKhodjaPrange12}
H.~A. Djikpesse, M.~R. Khodja, M.~D. Prange, S.~Duchenne, and H.~Menkiti.
\newblock Bayesian survey design to optimize resolution in waveform inversion.
\newblock {\em Geophysics}, 77(2):R81--R93, 2012.

\bibitem{EinsiedlerWard17}
M.~Einsiedler and T.~Ward.
\newblock {\em Functional analysis, spectral theory, and applications}, volume
  276 of {\em Graduate Texts in Mathematics}.
\newblock Springer, Cham, 2017.

\bibitem{EnglHankeNeubauer96}
H.~W. Engl, M.~Hanke, and A.~Neubauer.
\newblock {\em Regularization of {I}nverse {P}roblems}.
\newblock Springer Netherlands, 1996.

\bibitem{EtlingHerzog18}
T.~Etling and R.~Herzog.
\newblock Optimum experimental design by shape optimization of specimens in
  linear elasticity.
\newblock {\em SIAM Journal on Applied Mathematics}, 78(3):1553--1576, 2018.

\bibitem{FedorovLeonov13}
V.~V. Fedorov and S.~L. Leonov.
\newblock {\em Optimal design for nonlinear response models}.
\newblock CRC Press, 2013.

\bibitem{FengMarzouk19}
C.~Feng and Y.~M. Marzouk.
\newblock A layered multiple importance sampling scheme for focused optimal
  {B}ayesian experimental design.
\newblock {\em Preprint}, 2019.
\newblock \url{https://arxiv.org/abs/1903.11187}.

\bibitem{Fitzpatrick91}
B.~G. Fitzpatrick.
\newblock Bayesian analysis in inverse problems.
\newblock {\em Inverse problems}, 7(5):675, 1991.

\bibitem{Franklin70}
J.~N. Franklin.
\newblock Well-posed stochastic extensions of ill--posed linear problems.
\newblock {\em Journal of Mathematical Analysis and Applications}, 31:682--716,
  1970.

\bibitem{HaberHoreshTenorio08}
E.~Haber, L.~Horesh, and L.~Tenorio.
\newblock Numerical methods for experimental design of large-scale linear
  ill-posed inverse problems.
\newblock {\em Inverse Problems}, 24(055012):125--137, 2008.

\bibitem{HaberHoreshTenorio10}
E.~Haber, L.~Horesh, and L.~Tenorio.
\newblock Numerical methods for the design of large-scale nonlinear discrete
  ill-posed inverse problems.
\newblock {\em Inverse Problems}, 26(2):025002, 2010.

\bibitem{HaberMagnantLuceroEtAl12}
E.~Haber, Z.~Magnant, C.~Lucero, and L.~Tenorio.
\newblock Numerical methods for {A}-optimal designs with a sparsity constraint
  for ill-posed inverse problems.
\newblock {\em Computational Optimization and Applications}, pages 1--22, 2012.

\bibitem{HairerStuartVollmer14}
M.~Hairer, A.~M. Stuart, and S.~J. Vollmer.
\newblock Spectral gaps for a metropolis--hastings algorithm in infinite
  dimensions.
\newblock {\em The Annals of Applied Probability}, 24(6):2455--2490, 2014.

\bibitem{HannukainenHyvonenPerkkio20}
A.~Hannukainen, N.~Hyv{\"o}nen, and L.~Perkki{\"o}.
\newblock Inverse heat source problem and experimental design for determining
  iron loss distribution.
\newblock {\em Preprint}, 2020.
\newblock \url{https://arxiv.org/abs/2003.10395}.

\bibitem{HelinBurger15}
T.~Helin and M.~Burger.
\newblock Maximum a posteriori probability estimates in infinite-dimensional
  {B}ayesian inverse problems.
\newblock {\em Inverse Problems}, 31(8):085009, 2015.

\bibitem{Herman20}
E.~Herman.
\newblock {\em Design of inverse problems and reduced order modeling in complex
  physical systems}.
\newblock {PhD} dissertation, North Carolina State University, 2020.

\bibitem{HermanAlexanderianSaibaba19}
E.~Herman, A.~Alexanderian, and A.~K. Saibaba.
\newblock Randomization and reweighted $\ell_1$-minimization for {A}-optimal
  design of linear inverse problems.
\newblock {\em SIAM Journal on Scientific Computing}, 42(3):A1714--A1740, 2020.

\bibitem{HerzogRiedelUcinski18}
R.~Herzog, I.~Riedel, and D.~Uci{\'n}ski.
\newblock Optimal sensor placement for joint parameter and state estimation
  problems in large-scale dynamical systems with applications to
  thermo-mechanics.
\newblock {\em Optimization and Engineering}, 19(3):591--627, 2018.

\bibitem{HoangQuekSchwab19}
V.~H. Hoang, J.~H. Quek, and C.~Schwab.
\newblock Analysis of multilevel {MCMC-FEM} for {B}ayesian inversion of
  log-normal diffusions.
\newblock {\em SAM Research Report}, 2019, 2019.

\bibitem{HoreshHaberTenorio10}
L.~Horesh, E.~Haber, and L.~Tenorio.
\newblock Optimal experimental design for the large-scale nonlinear ill-posed
  problem of impedance imaging.
\newblock In {\em Large‐Scale Inverse Problems and Quantification of
  Uncertainty}, chapter~13, pages 273--290. John Wiley \& Sons, Ltd, 2010.

\bibitem{HuanMarzouk13}
X.~Huan and Y.~M. Marzouk.
\newblock Simulation-based optimal {B}ayesian experimental design for nonlinear
  systems.
\newblock {\em Journal of Computational Physics}, 232(1):288--317, 2013.

\bibitem{HuanMarzouk14}
X.~Huan and Y.~M. Marzouk.
\newblock Gradient-based stochastic optimization methods in {B}ayesian
  experimental design.
\newblock {\em International Journal for Uncertainty Quantification},
  4(6):479--510, 2014.

\bibitem{Hutchinson90}
M.~F. Hutchinson.
\newblock A stochastic estimator of the trace of the influence matrix for
  {L}aplacian smoothing splines.
\newblock {\em Communications in Statistics-Simulation and Computation},
  19(2):433--450, 1990.

\bibitem{HyvonenSeppanenStaboulis14}
N.~Hyvonen, A.~Sepp{\"a}nen, and S.~Staboulis.
\newblock Optimizing electrode positions in electrical impedance tomography.
\newblock {\em SIAM Journal on Applied Mathematics}, 74(6):1831--1851, 2014.

\bibitem{IsaacPetraStadlerEtAl15}
T.~Isaac, N.~Petra, G.~Stadler, and O.~Ghattas.
\newblock Scalable and efficient algorithms for the propagation of uncertainty
  from data through inference to prediction for large-scale problems, with
  application to flow of the {A}ntarctic ice sheet.
\newblock {\em Journal of Computational Physics}, 296:348--368, 2015.

\bibitem{JagalurMarzouk20}
J.~Jagalur-Mohan and Y.~Marzouk.
\newblock Batch greedy maximization of non-submodular functions: {G}uarantees
  and applications to experimental design.
\newblock 2020.
\newblock Preprint.~\url{https://arxiv.org/abs/2006.04554}.

\bibitem{KaipioSomersalo06}
J.~Kaipio and E.~Somersalo.
\newblock {\em Statistical and computational inverse problems}, volume 160.
\newblock Springer Science \& Business Media, 2006.

\bibitem{KieferWolfowitz59}
J.~Kiefer and J.~Wolfowitz.
\newblock Optimum designs in regression problems.
\newblock {\em The Annals of Mathematical Statistics}, 30:271--294, 1959.

\bibitem{KorkelBauerBockEtAl99}
S.~K{\"o}rkel, I.~Bauer, H.~G. Bock, and J.~Schl{\"o}der.
\newblock A sequential approach for nonlinear optimum experimental design in
  {DAE} systems.
\newblock {\em Scientific Computing in Chemical Engineering II}, 2:338--345,
  1999.

\bibitem{KorkelKostinaBockEtAl04}
S.~K{\"o}rkel, E.~Kostina, H.~G. Bock, and J.~P. Schl{\"o}der.
\newblock Numerical methods for optimal control problems in design of robust
  optimal experiments for nonlinear dynamic processes.
\newblock {\em Optimization Methods \& Software}, 19(3-4):327--338, 2004.
\newblock The First International Conference on Optimization Methods and
  Software. Part II.

\bibitem{KovalAlexanderianStadler19}
K.~Koval, A.~Alexanderian, and G.~Stadler.
\newblock Optimal experimental design under irreducible uncertainty for linear
  inverse problems governed by {PDEs}.
\newblock {\em Inverse Problems}, 36(7), 2020.

\bibitem{KrauseSinghGuestrin08}
A.~Krause, A.~Singh, and C.~Guestrin.
\newblock Near-optimal sensor placements in {G}aussian processes: Theory,
  efficient algorithms and empirical studies.
\newblock {\em Journal of Machine Learning Research}, 9:235--284, 2008.

\bibitem{Kullback1951}
S.~Kullback and R.~A. Leibler.
\newblock On information and sufficiency.
\newblock {\em The Annals of Mathematical Statistics}, 22(1):79--86, 03 1951.

\bibitem{LassasSaksmanSiltanen09}
M.~Lassas, E.~Saksman, and S.~Siltanen.
\newblock Discretization invariant {B}ayesian inversion and {B}esov space
  priors.
\newblock {\em Inverse Problems and Imaging}, 3(1):87--122, 2009.

\bibitem{LassasSiltanen04}
M.~Lassas and S.~Siltanen.
\newblock Can one use total variation prior for edge-preserving bayesian
  inversion?
\newblock {\em Inverse Problems}, 20:1537--1563, October 2004.

\bibitem{LehtinenPaivarintaSomersalo89}
M.~S. Lehtinen, L.~P\"{a}iv\"{a}rinta, and E.~Somersalo.
\newblock Linear inverse problems for generalized random variables.
\newblock {\em Inverse Problems}, 5:599--612, 1989.

\bibitem{Li19}
F.~Li et~al.
\newblock A combinatorial approach to goal-oriented optimal {B}ayesian
  experimental design.
\newblock Master's thesis, Massachusetts Institute of Technology, 2019.

\bibitem{LiuChepuriFardadEtAl16}
S.~Liu, S.~P. Chepuri, M.~Fardad, E.~Ma{\c{s}}azade, G.~Leus, and P.~K.
  Varshney.
\newblock Sensor selection for estimation with correlated measurement noise.
\newblock {\em IEEE Transactions on Signal Processing}, 64(13):3509--3522,
  2016.

\bibitem{LongMotamedTempone15}
Q.~Long, M.~Motamed, and R.~Tempone.
\newblock Fast {B}ayesian optimal experimental design for seismic source
  inversion.
\newblock {\em Computer Methods in Applied Mechanics and Engineering},
  291:123--145, 2015.

\bibitem{LongScavinoTemponeEtAl13}
Q.~Long, M.~Scavino, R.~Tempone, and S.~Wang.
\newblock Fast estimation of expected information gains for {B}ayesian
  experimental designs based on {L}aplace approximations.
\newblock {\em Computer Methods in Applied Mechanics and Engineering},
  259:24--39, 2013.

\bibitem{Mandelbaum84}
A.~Mandelbaum.
\newblock Linear estimators and measurable linear transformations on a
  {Hilbert} space.
\newblock {\em Zeitschrift f{\"u}r Wahrscheinlichkeitstheorie und Verwandte
  Gebiete}, 65(3):385--397, 1984.

\bibitem{McCormack18}
K.~A. McCormack et~al.
\newblock {\em Earthquakes, groundwater and surface deformation: exploring the
  poroelastic response to megathrust earthquakes}.
\newblock {PhD} dissertation, The University of Texas at Austin, 2018.

\bibitem{Muller05}
P.~M{\"u}ller.
\newblock Simulation based optimal design.
\newblock {\em Handbook of Statistics}, 25:509--518, 2005.

\bibitem{Muller07}
W.~G. M{\"u}ller.
\newblock {\em Collecting spatial data: optimum design of experiments for
  random fields}.
\newblock Springer Science \& Business Media, 2007.

\bibitem{NeitzelPieperVexlerEtAl19}
I.~Neitzel, K.~Pieper, B.~Vexler, and D.~Walter.
\newblock A sparse control approach to optimal sensor placement in
  {PDE}-constrained parameter estimation problems.
\newblock {\em To appear in Numerische Mathematik}, 2019.
\newblock \url{https://arxiv.org/abs/1905.01696}.

\bibitem{NemhauserWolseyFisher78}
G.~L. Nemhauser, L.~A. Wolsey, and M.~L. Fisher.
\newblock An analysis of approximations for maximizing submodular set
  functions—i.
\newblock {\em Mathematical programming}, 14(1):265--294, 1978.

\bibitem{Pazman86}
A.~P{\'a}zman.
\newblock {\em Foundations of Optimum Experimental Design}.
\newblock D. Reidel Publishing Co., 1986.

\bibitem{PetraMartinStadlerEtAl14}
N.~Petra, J.~Martin, G.~Stadler, and O.~Ghattas.
\newblock A computational framework for infinite-dimensional {B}ayesian inverse
  problems: {P}art {II}. {S}tochastic {N}ewton {MCMC} with application to ice
  sheet inverse problems.
\newblock {\em SIAM Journal on Scientific Computing}, 36(4):A1525--A1555, 2014.

\bibitem{PinskiSimpsonStuartEtAl15}
F.~J. Pinski, G.~Simpson, A.~M. Stuart, and H.~Weber.
\newblock Kullback--leibler approximation for probability measures on infinite
  dimensional spaces.
\newblock {\em SIAM Journal on Mathematical Analysis}, 47(6):4091--4122, 2015.

\bibitem{PratoZabczyk92}
G.~D. Prato and J.~Zabczyk.
\newblock {\em Stochastic Equations in Infinite Dimensions}.
\newblock Cambidge University Press, 1992.

\bibitem{PrenterVogel85}
P.~Prenter and C.~Vogel.
\newblock Stochastic inversion of linear first kind integral equations. {I}:
  Continuous theory and the stochastic generalized inverse.
\newblock {\em Journal of mathematical analysis and applications},
  106(1):202--218, 1985.

\bibitem{PronzatoPazman13}
L.~Pronzato and A.~P{\'a}zmanj.
\newblock {\em Design of Experiments in Nonlinear Models. Asymptotic Normality,
  Optimality Criteria amd Small-Sample Properties}.
\newblock New York: Springer-Verlag, 2013.

\bibitem{Pukelsheim93}
F.~Pukelsheim.
\newblock {\em Optimal Design of Experiments}.
\newblock John Wiley \& Sons, New-York, 1993.

\bibitem{ReedSimon72}
M.~Reed and B.~Simon.
\newblock {\em Methods of modern mathematical physics. {I}. {F}unctional
  analysis}.
\newblock Academic Press, New York-London, 1972.

\bibitem{RobertsRosenthal01}
G.~O. Roberts and J.~S. Rosenthal.
\newblock Optimal scaling for various metropolis-hastings algorithms.
\newblock {\em Statistical Science}, 16(4):pp. 351--367, 2001.

\bibitem{RudolfSprungk18}
D.~Rudolf and B.~Sprungk.
\newblock On a generalization of the preconditioned {Crank--Nicolson}
  metropolis algorithm.
\newblock {\em Foundations of Computational Mathematics}, 18(2):309--343, 2018.

\bibitem{ruthotto2017optimal}
L.~Ruthotto, J.~Chung, and M.~Chung.
\newblock Optimal experimental design for inverse problems with state
  constraints.
\newblock {\em SIAM Journal on Scientific Computing}, 40(4):B1080--B1100, 2018.

\bibitem{RyanDrovandiMcGreeEtAl16}
E.~G. Ryan, C.~C. Drovandi, J.~M. McGree, and A.~N. Pettitt.
\newblock A review of modern computational algorithms for {B}ayesian optimal
  design.
\newblock {\em International Statistical Review}, 84(1):128--154, 2016.

\bibitem{Ryan03}
K.~J. Ryan.
\newblock Estimating expected information gains for experimental designs with
  application to the random fatigue-limit model.
\newblock {\em Journal of Computational and Graphical Statistics},
  12(3):585--603, 2003.

\bibitem{SaibabaAlexanderianIpsen17}
A.~K. Saibaba, A.~Alexanderian, and I.~C. Ipsen.
\newblock Randomized matrix-free trace and log-determinant estimators.
\newblock {\em Numerische Mathematik}, 137(2):353--395, 2017.

\bibitem{SaibabaBardsleyBrownAlexanderian19}
A.~K. Saibaba, J.~Bardsley, D.~A. Brown, and A.~Alexanderian.
\newblock Efficient marginalization-based {MCMC} methods for hierarchical
  bayesian inverse problems.
\newblock {\em SIAM/ASA Journal on Uncertainty Quantification},
  7(3):1105--1131, 2019.

\bibitem{ShulkindHoreshAvron18}
G.~Shulkind, L.~Horesh, and H.~Avron.
\newblock Experimental design for nonparametric correction of misspecified
  dynamical models.
\newblock {\em SIAM/ASA Journal on Uncertainty Quantification}, 6(2):880--906,
  2018.

\bibitem{Sprungk20}
B.~Sprungk.
\newblock On the local {L}ipschitz stability of {B}ayesian inverse problems.
\newblock {\em Inverse Problems}, 36(5):055015, 2020.

\bibitem{Stuart10}
A.~M. Stuart.
\newblock Inverse problems: {A B}ayesian perspective.
\newblock {\em Acta Numerica}, 19:451--559, 2010.

\bibitem{SunseriHartVanBloemenWaandersEtAl20}
I.~Sunseri, J.~Hart, B.~van Bloemen~Waanders, and A.~Alexanderian.
\newblock {Hyper-Differential Sensitivity Analysis for Inverse Problems
  Constrained by Partial Differential Equations}.
\newblock {\em Inverse Problems}, 36(12):125001, 2020.

\bibitem{Tarantola05}
A.~Tarantola.
\newblock {\em Inverse Problem Theory and Methods for Model Parameter
  Estimation}.
\newblock SIAM, Philadelphia, PA, 2005.

\bibitem{Ucinski05}
D.~Uci{\'n}ski.
\newblock {\em Optimal measurement methods for distributed parameter system
  identification}.
\newblock CRC Press, Boca Raton, 2005.

\bibitem{Ucinski20}
D.~Uci{\'n}ski.
\newblock D-optimal sensor selection in the presence of correlated measurement
  noise.
\newblock {\em Measurement}, 2020.

\bibitem{Vogel02}
C.~R. Vogel.
\newblock {\em {C}omputational {M}ethods for {I}nverse {P}roblems}.
\newblock Frontiers in Applied Mathematics. Society for Industrial and Applied
  Mathematics (SIAM), Philadelphia, PA, 2002.

\bibitem{WalshWildeyJakeman17}
S.~Walsh, T.~Wildey, and J.~D. Jakeman.
\newblock Optimal experimental design using a consistent {B}ayesian approach.
\newblock {\em ASCE-ASME Journal of Risk and Uncertainty in Engineering
  Systems, Part B: Mechanical Engineering}, 2017.

\bibitem{Walter19}
D.~Walter.
\newblock {\em On sparse sensor placement for parameter identification problems
  with partial differential equations}.
\newblock PhD thesis, Technische Universit{\"a}t M{\"u}nchen, 2019.

\bibitem{WeiserFreytagErdmannEtAl18}
M.~Weiser, Y.~Freytag, B.~Erdmann, M.~Hubig, and G.~Mall.
\newblock Optimal design of experiments for estimating the time of death in
  forensic medicine.
\newblock {\em Inverse Problems}, 34(12):125005, 2018.

\bibitem{Williams1991}
D.~Williams.
\newblock {\em Probability with Martingales}.
\newblock Cambridge University Press, 1991.

\bibitem{WuChenGhattas20}
K.~Wu, P.~Chen, and O.~Ghattas.
\newblock A fast and scalable computational framework for large-scale and
  high-dimensional {B}ayesian optimal experimental design.
\newblock {\em Preprint}, 2020.
\newblock \url{https://arxiv.org/abs/2010.15196}.

\bibitem{YuZavalaAnitescu17}
J.~Yu, V.~M. Zavala, and M.~Anitescu.
\newblock A scalable design of experiments framework for optimal sensor
  placement.
\newblock {\em Journal of Process Control}, 2017.

\end{thebibliography}

\end{document}